\documentclass[12pt]{article}
\usepackage{amsmath,amsfonts,amssymb,amsthm,amscd}
\title{Groups, measures, and the NIP}
\author{Ehud Hrushovski\thanks{Supported by the Israel Science
Foundation grant no. 244/03}\\Hebrew University of Jerusalem
\and Ya'acov Peterzil\\University
of Haifa \and Anand
Pillay\thanks{Supported by NSF grants DMS-0300639
and the FRG DMS-0100979, as well as a Marie Curie chair}\\University 
of Illinois and University of Leeds}
\newtheorem{Theorem}{Theorem}[section]
\newtheorem{Proposition}[Theorem]{Proposition}
\newtheorem{Definition}[Theorem]{Definition}
\newtheorem{Remark}[Theorem]{Remark}
\newtheorem{Lemma}[Theorem]{Lemma}
\newtheorem{Corollary}[Theorem]{Corollary}

\newtheorem{Example}[Theorem]{Example}
\newtheorem{Question}[Theorem]{Question}
\newcommand{\R}{\mathbb R}   \newcommand{\Rr}{\mathbb R}
\newcommand{\Q}{\mathbb Q}  \newcommand{\Qq}{\mathbb Q}
\newcommand{\Z}{\mathbb Z}  
\newcommand{\N}{\mathbb N}  \newcommand{\Nn}{\mathbb N}
\newcommand{\Pp}{\mathbb P}
\def\union{\cup}
\newcommand{\pf}{\noindent{\em Proof. }}
\newcommand{\la}{\langle}
\newcommand{\ra}{\rangle}
\begin{document}
\maketitle

\begin{abstract} We discuss measures, invariant measures on
definable groups, and genericity, often in an NIP (failure of the
independence property) environment. We complete the proof of the
third author's conjectures relating definably compact groups $G$ in
saturated $o$-minimal structures to compact Lie groups. We also
prove some other structural results about such $G$, for example the
existence of a left invariant finitely additive probability measure
on definable subsets of $G$. We finally introduce a new notion
``compact domination" (domination of a definable set by a compact
space) and raise some new conjectures in the $o$-minimal case.
\end{abstract}

\section{Introduction}
\label{1} One of the occasions for writing this paper is the
completion of the proof of the ``$o$-minimal group conjectures" of
the third author, from \cite{Pillay}. Among the new ingredients are
(i) the use of invariant measures on definable sets in the presence
of  the NIP (failure of the independence property), and (ii) the
identification of a certain property (finitely satisfiable generics)
which can be used in an inductive proof, and is of interest in its
own right.

The measures appear in Keisler's paper \cite{Keisler} which
is a strong influence on our work. In Keisler's work, the
theory of forking is in a sense extended from stable theories
to theories without the independence property, but replacing
complete types by measures (on the Boolean algebra of
definable sets). It is somewhat amusing to note that
Keisler's work was roughly contemporaneous with early work
on
$o$-minimality which was also motivated by the attempt to
generalize stability to suitable ordered structures.

Our work may also overlap to some extent with recent papers
of Shelah on theories without the independence property
(for example \cite{Shelah1}, \cite{Shelah}).

In any case, we take the opportunity in this paper to expand on and
develop some theory, not all of which is directed towards the proof
of the $o$-minimal group conjectures.

Stability and stable group theory are at the core of ``pure" or
``abstract" model theory. Recall Shelah's result that $T$ is stable
iff $T$ does not have the strict order property and does not have
the independence property (see \cite{Shelah1}. There has been
considerable work on generalizing stability to particularly nice
theories without the strict order property, namely the simple
theories. So part of this paper is around developing some theory in
an ``orthogonal'' direction, namely for certain theories $T$ without
the independence property. Another aspect of this paper is the
``model theory of the standard part map".

In Section 2, we recall and elaborate on some of Keisler's notions
from \cite{Keisler}. In particular we discuss smooth, definable, and
finitely satisfiable measures. In Section 3, we discuss some
consequences of NIP, sometimes in the presence of measures. Include
here is a ``Borel definability" of coheirs assuming NIP. In Section
4, we introduce the ``finitely satisfiable generics" property for
definable groups $G$, stating which aspects of stable group theory
are valid in this situation. In Section 5 we discuss in general
``definably amenable groups'', namely groups with a left invariant
measure on the definable sets. In Section 6 we prove various results
around existence of $G^{00}$ and existence of invariant measures
under the NIP assumption. In Section 7 we take a short diversion to
explain how our results can generalize to the class of ``inductively
definable" groups. In Section 8 we prove the full conjecture from
\cite{Pillay}:

(*) {\em If $G$ is a definably compact group definable in a saturated
$o$-minimal expansion of a real closed field, then the quotient
$G/G^{00}$ of $G$ by its smallest type-definable subgroup of bounded
index $G^{00}$, is, when equipped with the logic topology,  
a 
compact Lie group whose dimension (as a Lie
group) equals the dimension of $G$ (as a definable set in an
$o$-minimal structure).}

  The proof rests on and continues a number of earlier papers
\cite{Pillay}, \cite{BO}, \cite{BOPP}, \cite{Peterzil-Pillay}, and
\cite{Edmundo-Otero}. We will give below a guide
for the reader who 
is interested in a fast path to the proof of (*).

In Section 9 and 10, we isolate a new notion,  of ``compact
domination", and conjecture that in fact a definably compact group
$G$ in an $o$-minimal structure is compactly dominated by
$G/G^{00}$. We then prove this in several special 
cases.

\vspace{5mm}
\noindent
{\em Guide to the proof of (*).}  The 
proof is carried out in section 8. Globally it proceeds by induction 
on $dim(G)$. The two extreme cases are when (a) $G$ is commutative, 
and (b) $G$ is definably simple.
The ``new" ingredient for case (a) 
is use of the amenability of $G$ (namely the existence of an 
invariant finitely additive measure on all subsets of $G$) together 
with the $NIP$. The key sequence of preliminary results is Lemma 2.8, 
Proposition 3.3, Corollary 3.4 and Proposition 6.3. Case (a) is 
proved in Lemma 8.2. Case (b) was proved in \cite{Peterzil-Pillay} 
under the weaker hypothesis of  ``$G$ has very good reduction". This 
is discussed in Lemma 8.3 of the current paper. For the induction 
step, one may assume $G$ has a normal commutative definable subgroup 
$N$. But we need to know more than simply that (*) holds for $G/N$ 
and $N$. Namely we require that both $G/N$ and $N$ have the 
``finitely satisfiable generics" property. The $fsg$ is introduced in 
section 4, and Proposition 4.2 is crucial.  In Cases (a) and (b) we 
actually prove in addition that the relevant groups have the $fsg$ 
property. Proposition 4.5 shows that from the $fsg$ for $G/N$ and $N$ 
we can conclude the $fsg$ for $G$. An argument using Corollary 4.3 
shows that (*) holds for $G$.

\vspace{5mm} \noindent 
Our notation is standard. We work in a large
saturated model ${\bar M}$ of a complete first order, possibly
many-sorted theory $T$ in a language $L$. If we assume that  $|{\bar
M}| = {\bar\kappa}$ then by a ``small" or ``bounded" set we mean a set
of cardinality $< {\bar \kappa}$. $x, y$ denote finite sequences of
variables unless we say otherwise. $A,B,..$ denote small subsets of
${\bar M}$. $M, N,..$ denote small elementary substructures of
${\bar M}$. ``Type-definable" means the intersection of a small
collection of definable sets, and a ``bounded type-definable
equivalence relation'' is a type-definable equivalence relation with
a bounded number of classes. We refer to \cite{Pillay-book} for any
background on stability.

$T$ is said to have the NIP (for ``not the independence property")
if there is no formula $\phi(x,y)\in L$ and $\la a_{i}:i<\omega\ra$
and $\la b_{w}:w\subseteq \omega\ra $ such that $\models
\phi(a_{i},b_{w})$ iff $i\in \omega$. Stable and o-minimal theories,
as well as the theory of the p-adic field are all examples of
theories with NIP, while simple unstable theories all have the
independence property.

If $G$ is a group definable in ${\bar M}$ then $G^{00}$ is the
smallest type-definable subgroup of bounded index in $G$, if there
is such. If $E$ is a type-definable equivalence relation on a
definable set $X$ with a bounded number of classes, then the logic
topology on $X/E$ is given by: $C\subseteq X/E$ is closed if the
pre-image of $C$ in $X$ is type-definable.

In various parts of the 
paper we will make use of standard facts and techniques
regarding 
indiscernibles, which the referee has asked us to explain.  One of 
these facts is that given a complete theory $T$, and cardinal $\mu$ 
there is a cardinal $\lambda$ such that if 
$\{a_{\alpha}:\alpha<\lambda\}$ is a set of $\mu$-tuples in some 
saturated model of $T$, then there is an indiscernible sequence 
$(b_{i}:i<\omega)$ of $\mu$-tuples, such that for every $n$, 
$tp(b_{0},..,b_{n-1}) = tp(a_{\alpha_{0}},....,a_{\alpha_{n-1}})$ for 
some $\alpha_{0} < .. < \alpha_{n-1} < \lambda$. This is an 
application of the Erd\"os-Rado Theorem. A statement and proof 
appears in \cite{Grossberg} (Theorem 1.13) for example. When using 
this fact we will just say ``by Erd\"os-Rado". Another method is 
``stretching" indiscernibles: namely given an indiscernible sequence 
$(a_{i}:i<\omega)$ we can, for any totally ordered set $I$, find an 
indiscernible sequence $(b_{i}:i\in I)$ such that for each $n$ and 
$i_{0} < ...< i_{n}$ in $I$, $tp(b_{i_{0}},...,b_{i_{n-1}}) = 
tp(a_{0},..,a_{n-1})$. This is of course just by compactness.

\vspace{5mm}
\noindent
Some of the work presented here was done while the authors
were at the Isaac Newton Institute, Cambridge, for the
Spring 2005 Model Theory program. We would like to thank both the
Newton Institute and the organizers of the program for their
hospitality, ideal conditions and financial support. In addition to 
the referee, several other individuals
and research groups have 
passed on to us comments on an earlier version of the paper as well 
as helpful suggestions. So we would also like
to thank Alessandro 
Berarducci, Margarita Otero, and Lou van den Dries and
participants 
in the UIUC model theory seminar.

\section{Definable functions and measures}
\label{2} We consider here functions of one kind or another from
sorts, or definable sets in ${\bar M}$, to compact Hausdorff spaces
$C$, such as the closed interval $[0,1]$.
\begin{Definition} \label{2.1} Let $X$ be an $A$-definable set in $M$, $C$ some
compact Hausdorff space of bounded size, and $f$ a map from $X$ to 
$C$. We will say that
$f$ is {\em definable over $A$}, if for any closed subset $C_{1}$ of $C$,
$f^{-1}(C_{1}) \subseteq X$ is type-definable over $A$ in $M$.
\end{Definition}

\begin{Example}\label{2.2} (i) The {\em tautological} map $s$ from
$X$ to its Stone space $S_{X}(A)$: $s(b) = tp(b/A)$. Note that a map
$f$ from
$X$ to a compact Hausdorff space $C$ will be definable over
$A$ just if $f = g\circ s$ with $g$ a continuous map from
$S_{X}(A)$ to $C$. So the tautological definable map $s$ is
also universal.
\newline
(ii) Let $A$ be a small subset of sort $X$ in ${\bar M}$, and
$\phi(x,y)$ a formula, with $x$ of sort $X$ and $y$ of sort $Y$.
Identify the power set of $A$ with the compact space $2^{|A|}$. Let
$f:Y\rightarrow 2^{|A|}$ be given by $f(b) = \{a\in A:\models
\phi(a,b)\}$. Then, as is easy to verify,  $f$ is definable over
$A$.
\end{Example}

In Definition \ref{2.1}, note that if $f:X\rightarrow C$ is
definable, then $f(X)\subseteq C$ is closed  (because as in Example
\ref{2.2}(i), $f$ can be identified with a continuous map between
compact spaces hence its image is closed).  So we may assume $f$ to
be onto.

In fact definable maps as in Definition \ref{2.1} amount to the same
thing as quotienting by bounded type-definable equivalence
relations:
\begin{Remark} \label{2.3} Let $X$ be definable over $A$ in $M$.
\newline
(i) Let $f$ be a definable (over $A$) map from $X$ onto the compact
Hausdorff space $C$ in the sense of Definition \ref{2.1}. Let $E =
\{(x,y)\in X\times X: f(x) = f(y)\}$. Then $E$ is an
$A$-type-definable equivalence relation of bounded index, and $f$
induces a homeomorphism between $X/E$ with the logic topology and
the space $C$.
\newline
(ii) Conversely, if $E$ is a bounded $A$-type-definable equivalence
relation on $X$, $X/E$ is equipped with the logic topology, and
$M_{0}$ is a small model containing $A$ and a representative for
each $E$-class, then the quotient map $f:X\rightarrow X/E$ is an
$M_{0}$-definable map from $X$ onto the compact Hausdorff space
$X/E$.
\end{Remark}
\pf (i) For each pair $C_{1}, C_{2}$ of closed subsets of $C$ such
that $C_{1}\cup C_{2} = C$, let $E_{C_{1}, C_{2}} = \{(x,y)\in
X\times X$ such that either $f(x)\in C_{1}$ and $f(y)\in C_{1}$ or
$f(x)\in C_{2}$ and $f(y)\in C_{2}$\}. So $E_{C_{1},C_{2}}$ is
type-definable over $A$. As $X$ is Hausdorff, $E$ is the
intersection of all $E_{C_{1},C_{2}}$ hence is also type-definable.
Identifying $X/E$ with $C$ we see that the logic topology on $C$
refines the original topology on $C$. As both topologies are compact
Hausdorff they agree. $E$ is of bounded index since the pre-image of
each singleton in $C$ is type-definable over a fixed set $A$.
\newline
(ii) If $C\subseteq X/E$ is closed, then by definition $f^{-1}(C)$ is 
type-definable.
But $f^{-1}(C)$ is also $M_{0}$-invariant, hence it 
is type-definable over $M_{0}$. 
\qed

\vspace{5mm} \noindent So Definition \ref{2.1} is cosmetic. However
it enables some unification of various notions, as well as some
clean statements. For example the conjecture from \cite{Pillay} can
now be restated as:

{\em If $G$ is a definably connected definably compact group in a
saturated $o$-minimal structure $M$ then there is a definable
surjective homomorphism $f$ from $G$ to a compact Lie group $G_{1}$
where $dim(G_{1})$ equals the $o$-minimal dimension of $G$. Moreover
any other definable homomorphism from $G$ into a compact group
factors through $f$}.

\vspace{5mm}
\noindent
We now recall the probablity measures on definable sets considered by
Keisler \cite{Keisler}. We will call these {\em Keisler measures}. Let us
fix again a sort or definable set $X$ in ${\bar M}$ which we assume to be
$\emptyset$-definable. $Def(X)$ will denote the subsets of $X$ definable
(with parameters) in $\bar M$, and $Def_{A}(X)$ those sets defined over
$A$. (So we identify $Def(X)$ with $Def_{\bar M}(X)$.)
\begin{Definition}\label{2.4}{\em  (i) A {\em Keisler measure} $\mu$ 
on $X$ over $A$ is a
finitely additive probability measure on $Def_{A}(X)$; namely a map $\mu$
from $Def_{A}(X)$ to the interval $[0,1]$ such that $\mu(\emptyset) = 0$,
$\mu(X) = 1$ and for $Y,Z\in Def_{A}(X)$, $\mu(Y\cup Z) = \mu(Y) +
\mu(Z) - \mu(Y\cap Z)$.
\newline
(ii) {\em A (global) Keisler measure} on $X$ is a finitely additive
probability measure on $Def(X)$.
\newline
(iii) If $\mu$ is a Keisler measure on $Def_{B}(X)$ and $A\subseteq
B$ we write $\mu|A$ for the restriction of $\mu$ to $Def_{A}(X)$.}
\end{Definition}

  Note that a complete type (of an element of $X$) over $A$ is
precisely a $0$-$1$ valued Keisler measure on $X$ over $A$.

For each $L$-formula $\phi(x,y)$ with $x$ a variable of sort $X$,
let $S_{\phi}$ be the sort whose elements are the subsets of $X$
defined by instances of $\phi$. So a global Keisler measure on $X$
is given through a family $\{\mu_{\phi}:\phi(x,y)\in L\}$ of maps
$\mu_{\phi}:S_{\phi}\rightarrow [0,1]$.

Keisler observes that any Keisler measure on $X$ over $A$ extends to
a global Keisler measure on $X$. Moreover any Keisler measure on $X$
over $A$ extends to a unique countably additive measure on the
$\sigma$-algebra generated by the $A$-definable subsets of $X$ (see
Theorem 1.2 in \cite{Keisler}). We will point out now a way of
extending a Keisler measure over a model to a global Keisler
measure, as the construction will be useful later on. \vspace{3mm}

\noindent {\bf Construction }(*){\em  Let $\mu$ be a Keisler measure
on $X$ over a model $M_{0}$, viewed as a map from definable in
$M_{0}$ subsets of $X(M_{0})$ to $[0,1]$. Consider the structure
$\la M_{0},[0,1],+,<,\mu_{\phi}\ra_{\phi}$ consisting of
$M_{0}^{eq}$, the real unit interval $[0,1]$, and for each $\phi$,
the map $\mu_{\phi}:S_{\phi}(M_{0})\rightarrow [0,1]$ as well as the
ordering and addition (modulo $1$) on $[0,1]$. Take a saturated
elementary extension $\la M_{0}',[0,1]',+,<, \mu'_{\phi}\ra_{\phi}$.
Then the composition of $\mu'$ with the standard part map
$st:[0,1]'\rightarrow [0,1]$ is a Keisler measure on $X$ over
$M_{0}'$ extending $\mu$. We may identify ${\bar M}$ with $M_{0}'$}.

\vspace{2mm} \noindent One point of this construction is that the
structure $\bar M$, equipped with the constructed measure, has some
obvious ``saturation" properties.

\vspace{2mm} \noindent We have observed that a Keisler measure on
$X$ is (among other things) a sequence of maps from sorts $S_{\phi}$
to $[0,1]$. It would be natural to call $\mu$ {\em definable} if
each $\mu_{\phi}:S_{\phi}\rightarrow [0,1]$ is definable in the
sense of Definition \ref{2.1}. This is precisely (i) in the next
definition.

\begin{Definition}{\em  Let $\mu$ be a (global) Keisler measure on $X$.
\newline
(i) Then $\mu$ is {\em definable over $A$} iff for each $L$-formula
$\phi(x,y)$, and closed subset $C$ of $[0,1]$, $\{b\in
M:\mu(\phi(x,b))\in C\}$ is type-definable over $A$.

  Let $M_{0}$ be
a small submodel of ${\bar M}$.

\noindent (ii) We say that $\mu$ is {\em finitely satisfiable} in
$M_{0}$ if whenever $Y\subseteq X$ is definable and $\mu(Y)>0$ then
$Y\cap M_{0} \neq \emptyset$.

\noindent (iii) We say that $\mu$ is {\em smooth over $M_{0}$} if
$\mu$ is the unique (global) extension of $\mu|M_{0}$ to a measure
on $X$. In this situation we also say that $\mu|M_{0}$ is smooth.}
\end{Definition}

The notion of a smooth measure was also introduced by Keisler
(\cite{Keisler}) although his definition is weaker than the above,
for certain technical reasons. In any case, if $\mu$ is a $0-1$
measure given by a complete type then it is smooth if and only if
the type is algebraic.

Here is a ``nonalgebraic'' example of a smooth Keisler measure: Let
$\bar M$ be a saturated real closed field, and take $X$ to be the
interval $[0,1]$ in the sense of ${\bar M}$. The field of reals $\R$
is an elementary substructure of ${\bar M}$. The standard measure on
the real unit interval $[0,1]^{\R}$ gives a Keisler measure on $X$
over $\R$ which is easily seen to have a unique extension over
${\bar M}$. (This will be subsequently generalized in the last
section.)

\begin{Lemma}\label{2.6} Let $\mu$ be a (global) Keisler measure on 
$X$. Suppose
that $\mu$ is smooth over $M_{0}$. Then $\mu$ is both finitely
satisfiable in $M_{0}$ and definable over
$M_{0}$.
\end{Lemma}
\pf Finite satisfiability is immediate from \cite{Keisler}, Lemma
2.2 (which is itself based on Lemma 1.6 there), but for the sake of
completeness we repeat the argument here.

It is clearly sufficient to prove that if $X$ is a definable set in
$\bar M$ with $\mu(X)>0$ then it contains an $M_0$-definable $Y$
with $\mu(Y)>0$. Assume not, namely that all $M_0$-definable subsets
of $X$ have $\mu$-measure zero. By the smoothness assumption, it is
sufficient to show that there is {\em some} finitely additive
Keisler-measure $\mu'$ on $\bar M$, extending $\mu|M_0$, with
$\mu'(X)=0$. By compactness, this amounts to showing, given finitely
many $M_0$-definable sets $Y_1,\ldots,Y_k$, that there is a finitely
additive probability measure $\mu'$ on the Boolean algebra generated
by $Y_1,\ldots,Y_k,X$, which agrees with $\mu$ on the $Y_i$'s. Let
$\mathcal B_0$ be the Boolean algebra generated by the $Y_i$'s.
Without loss of generality, the $Y_i$'s are atoms in $\mathcal B_0$
and hence each $Y_i\cap X$ is an atom in the Boolean algebra
generated by $\mathcal B_0$ and $X$. We now let $\mu'(Y)=\mu(Y)$ for
all $Y\in \mathcal B_0$ and $\mu'(Y_i\cap X)=0$. This gives the
desired measure $\mu'$ and proves that $\mu$ is finitely
satisfiable.

The definability of $\mu$ over $M_{0}$ is more or less explained by
a ``Beth's Theorem for continuous logic". But we will be more
direct. We make use of Construction (*) above. Consider the
structure $\la M_{0},[0,1], +,<, \mu_{\phi}|M_{0}\ra_{\phi}$ from
there, equipped with constants for all elements (of $M_{0}$ and of
the unit interval). Let $T_{1}$ be its theory. We saw that in a
saturated model $\bar M_1$ of $T_{1}$, $\{st\circ\mu'_{\phi}:
\phi\in L\}$ gives rise to a Keisler measure $\mu''$ extending
$\mu|M_{0}$.  We may assume that $\bar M_1$ is an expansion of $\bar
M$, and by the smoothness assumption, that $\mu''=\mu$.

Fix an $L$-formula $\phi(x,y)$ where $x$ is of sort $X$. Given a
closed set $C\subseteq [0,1]$, we want to show that the set
$X_1=\{b: \mu(\phi(x,b))\in C\}$ is type-definable in $\bar M$ over
$M_0$. Note that the standard part map $st:[0,1]'\to [0,1]$ (where
$[0,1]'$ is the unit interval in $\bar M_1$) is definable in $\bar
M_1$ (over the empty set) in the sense of Definition \ref{2.1}, and
by the definability of $\mu$ in $\bar M_1$, the set $X_1$ is
type-definable over $M_0$ in $\bar M_1$, via a type $\Sigma(y)$.

Now, the smoothness assumption implies that $\Sigma(y)$ does not
depend on the particular expansion $\bar M_1$ of $\bar M$. We can
now apply the classical Beth  Theorem (for types rather than
formulas) and conclude that $X_1$ is type-definable in $\bar M$,
over $M_0$.\qed


\begin{Remark} Let $\mu$ be a global Keisler measure on $X$. Let us
define $\mu$ to be an heir of $\mu|M_{0}$ if for each $L$-formula
$\phi(x,y)$ and $r\in [0,1)$, if for some $b\in {\bar M}$,
$\mu(\phi(x,b))> r$ then for some $b\in M_{0}$,
$\mu(\phi(x,b))> r$. Then the proof above shows that $\mu$
is the unique heir of $\mu|M_{0}$ over $\bar M$ if and only
if $\mu$ is definable over $M_{0}$.
\end{Remark}

\vspace{2mm}
\noindent
The following relationship between Keisler measures and indiscernibles
will be useful. It also appears in
\cite{Newelski-Petrykowski}.
\begin{Lemma}\label{2.8} Let $\mu$ be a Keisler measure on $X$. Let $x$ be a
variable of sort $X$, let $\phi(x,y)\in L$, and let $\la
b_{i}:i<\omega\ra$ be an indiscernible sequence such that for some
$\epsilon > 0$, $\mu(\phi(x,b_{i})) \geq \epsilon$ for all $i$. Then
$\{\phi(x,b_{i}): i<\omega\}$ is consistent.
\end{Lemma}
\pf Let $Y_{b_{i}}$ denote the set defined by $\phi(x, b_{i})$. By
construction (*) above and Ramsey's theorem, we may assume that the
sequence $\la b_{i}:i<\omega\ra$ is also indiscernible with respect
to the map $\mu$, in particular that for each
$i_{1}<..<i_{n}<\omega$ and $j_{1} < ..< j_{n} < \omega$,
$\mu(Y_{b_{i_{1}}}\cap ...\cap Y_{b_{i_{n}}}) =
\mu(Y_{b_{j_{1}}}\cap ...\cap Y_{b_{j_{n}}}) = r_{n}$ say. So by
assumption $r_{1} > 0$.

Suppose for a contradiction that some finite intersection of the
$Y_{b_{i}}$'s is empty. Choose maximal $k$ such that $r_{k}> 0$. For
$j\geq 0$ let $Z_{j} = Y_{b_{1}}\cap Y_{b_{2}} \cap ..\cap
Y_{b_{k-1}}\cap Y_{b_{k+j}}$. Then each $Z_{j}$ has measure
$r_{k}>0$ and their pairwise intersections have measure $0$, a
contradiction.\qed

\section{NIP and some consequences}
The definition of NIP (failure of independence property) was given
in the Introduction. A well-known equivalence (see Theorem 12.17 of 
\cite{Poizat}) is:
\begin{Lemma} \label{3.1} $T$ has the NIP if and only for any sequence
$\la b_{i}:i<\omega\ra$ which is indiscernible over $\emptyset$ and
formula $\phi(y)$, possibly with parameters, there is an $i$ such
that $\models \phi(b_{j})$ for all $j>i$, or $\models\neg\phi(b_{j})$ 
for all $j<i$.
\end{Lemma}

Notation: if $\phi(x), \psi(x)$ are formulas, let
$\phi(x)\Delta\psi(x)$ denote the symmetric difference
$(\phi(x)\wedge\neg\psi(x))\vee(\neg\phi(x)\wedge\psi(x))$
of $\phi$ and $\psi$.

\begin{Corollary} \label{3.2} Suppose $T$ has NIP. Let $\phi(x,y)$ be an
$L$-formula, and $\la b_{i}:i<\omega\ra$ an indiscernible sequence.
Then the set  $\{\phi(x,b_{2j})\Delta\phi(x,b_{2j+1}):j<\omega\}$ is
inconsistent.
\end{Corollary}
\pf Otherwise, let $c$ realize
$\{\phi(x,b_{2j})\Delta\phi(x,b_{2j+1}):j<\omega\}$ and the formula
$\phi(c,y)$ contradicts Lemma \ref{3.1}.\qed

\vspace{5mm}
\noindent
We now give some consequences of the NIP for Keisler measures. The main
insight is due to Keisler (\cite{Keisler}, Theorem 3.14). We are back to
the context of
$\bar M$ a saturated model of $T$ and $X$ a sort or $\emptyset$-definable
set in $\bar M$.
\begin{Proposition} \label{3.3} Assume $T$ has the NIP. Let $\mu$ be a (global)
Keisler measure on $X$. Let $\phi(x,y)$ be a formula with $x$ of
sort $X$, and $\epsilon > 0$. Then there do not exist $\la
b_{i}:i<\omega\ra$ such that $i\neq j$ implies
$\mu(\phi(x,b_{i})\Delta\phi(x,b_{j})\geq \epsilon$.
\end{Proposition}
\pf Suppose otherwise. Then by Construction (*) from Section 2, and
Ramsey's theorem, we may assume in addition that $\la
b_{i}:i<\omega\ra$ is indiscernible. By Lemma \ref{2.8},
$\{\phi(x,b_{2j})\Delta\phi(x,b_{2j+1}):j< \omega\}$ is consistent,
contradicting Corollary \ref{3.2}.\qed

\begin{Corollary} \label{3.4}Assume $T$ has NIP and let $\mu$ be a 
global Keisler
measure on $X$. For definable subsets $Y, Z$ of $X$, define
$Y\sim_{\mu}Z$ if $\mu(Y\Delta Z) = 0$. Then there are only boundedly
many $\sim_{\mu}$-classes of definable subsets of $X$. In particular
there is a small model $M_{0}$ such that every definable subset $Y$ of
$X$ is $\sim_{\mu}$ to some $M_{0}$-definable subset of $X$.
\end{Corollary}
\pf If there are unboundedly many definable subsets of $X$ modulo
$\sim_{\mu}$ then we can clearly find a formula $\phi(x,y)$ and
large set $\la b_{i}:i\in I\ra$ such that the measures of the
pairwise symmetric differences of the $\phi(x,b_{i})$ are $>0$. By
Construction (*) from Section 2, we may assume that $\la b_{i}:i\in
I\ra$ is an indiscernible sequence with respect to $\mu$ as well,
whereby $\mu(\phi(x,b_{i})\Delta\phi(x,b_{j}))\geq \epsilon$ for
some fixed $\epsilon > 0$ and all $i\neq j$. This contradicts
Proposition \ref{3.3}.\qed

\vspace{5mm}
\noindent
Our next result is in a somewhat different spirit.
\begin{Theorem} \label{3.5} Suppose $T$ is countable with NIP. Let
$M_{0}$ be a countable elementary substructure of $\bar M$. Let
$p(x)$ be a complete $1$-type over ${\bar M}$ which is finitely
satisfiable in $M_{0}$. Let $U = \{X\cap M_{0}:X\in
p\}$. Then $U$ is a Borel (in fact an $F_{\sigma}$) subset of the 
Polish space
$2^{M_{0}}$.
\end{Theorem}

Before going into the proof we give an easy example to
illustrate the technique.
\begin{Remark} \label{3.6} Let $T$ be countable, and let $M_{0}$ be a
countable model. Then
\newline
(i) The set $\{X\cap M_{0}: X$ a definable subset of ${\bar
M}\}$ is an $F_{\sigma}$ (as a subset of $2^{M_{0}}$).
\newline
(ii) Let $p(x)\in S_{1}({\bar M})$ be definable. Then
$\{X\cap M_{0}: X\in p\}$ is an $F_{\sigma}$.
\end{Remark}
\pf (i)  Fix an
$L$-formula $\phi(x,y)$, and $n<\omega$. Let $U_{\phi} = \{X\cap
M_{0}: X$ is defined by $\phi(x,c)$ for some $c$\}. By Example
\ref{2.2}(ii), $U_{\phi}$ is closed. Then $U = \cup_{\phi}U_{\phi}$
is Borel and coincides with $\{X\cap M_{0}:X$ definable subset of
${\bar M}\}$.
\newline
(ii) Suppose again $\phi(x,y)\in L$ and let $\psi(y,d)$ be a formula
defining $p|\phi$. Then define $U_{\phi}$ just as above but
requiring also that $c$ realizes $\psi(y,d)$.\qed

\vspace{5mm} \noindent The proof of Theorem \ref{3.5} will go
through several lemmas.

For now let $T$ be an arbitrary complete theory with NIP.
\begin{Lemma}\label{3.7} For any $\phi(x,y)\in L$, there is some $N =
N_{\phi}$, such that for any indiscernible sequence $\la
a_{i}:i<\omega\ra$ and $c$, there do not exists $i_{0} < i_{1} < ..
< i_{N}$ such that for each $j<N$, $\models
\phi(a_{i_{j}},c)\leftrightarrow \neg\phi(a_{i_{j+1}},c)$.
\end{Lemma}
\pf Otherwise, by compactness we find an indiscernible sequence $\la
a_{i}:i<\omega\ra$ and $c$ such that for each $i<\omega$, $\models
\phi(a_{i},c)$ iff $\models \neg\phi(a_{i+1},c)$, contradicting
Lemma \ref{3.1}.\qed

\vspace{5mm} \noindent Recall that a type $p(x)\in S(\bar M)$ is
called {\em finitely satisfiable} in a model $M_0\subseteq \bar M$
if every formula in $p(x)$ is satisfiable in $M_0$.  If $p(x)\in
S(\bar M)$ is finitely satisfiable in a small model $M_{0}$, then we
can build an indiscernible sequence $I = \la a_{0},a_{1},...\ra$
over $M_{0}$ by letting $a_{0}$ realize $p|M_{0}$ and $a_{i+1}$
realize $p|(M_{0}\cup\{a_{0},..a_{i}\})$. Although the sequence $I$
is not unique, its type $tp(\la a_{i}:i<\omega\ra /M_{0})$ IS
unique, and we call this type $Q_{p,M_{0}}$.

Let us now fix a type $p(x)\in S(\bar M)$ which is finitely
satisfiable in $M_{0}$, and let $Q = Q_{p,M_{0}}$. (So $Q$
is a complete type over $M_{0}$ in variables
$(x_{i}:i<\omega)$). Let $Q_{n}$ be the restriction of $Q$
to the variables $(x_{0},..,x_{n})$. Fix an $L$-formula
$\phi(x,y)$ and some $c$ from $\bar M$. We will say that a
realization $(a_{0},..,a_{n})$ of $Q_{n}$ is {\em good} for
$\phi(x,c)$, if
\newline
(i) $\models
\phi(a_{i},c)\leftrightarrow\neg\phi(a_{i+1},c)$ for all
$i<n$, and
\newline
(ii) there does not exist $a_{n+1}$ such that
$(a_{0},..,a_{n},a_{n+1})$ realizes $Q_{n+1}$ and
\newline
$\models \phi(a_{n},c)\leftrightarrow \neg\phi(a_{n+1},c)$.

\vspace{2mm}
\noindent

With this notation, we have the following:
\begin{Lemma}\label{3.8} For $p$ as above, the following are equivalent:
\newline
(i) $\phi(x,c)\in p$,
\newline
(ii) there is $k\leq N_{\phi}$ and there is a realization
$(a_{0},..,a_{k})$ of $Q_{k}$ which is good for $\phi(x,c)$ such
that $\models \phi(a_{k},c)$,
\end{Lemma}
\pf Note first that by Lemma \ref{3.7}, for any $c$ there is $k\leq
N_{\phi}$ and realization $(a_{0},..,a_{k})$ of $Q_{k}$ which is
good for $\phi(x,c)$.
\newline
Now suppose $(a_{0},..,a_{k})$ realizes $Q_{k}$ and is good for
$\phi(x,c)$. Let $M_{1}$ be a small model containing
$M_{0}\cup\{a_{0},..a_{k},c\}$ and let $a$ realize $p|M_{1}$. Note
that $(a_{0},..,a_{k},a)$ realizes $Q_{k+1}$. By the ``goodness" of
$(a_{0},..,a_{k})$ for $\phi(x,c)$, it follows that $\models
\phi(a_{k},c)\leftrightarrow \phi(a,c)$. But $\models \phi(a,c)$ iff
$\phi(x,c)\in p$.

This is enough to prove the lemma.\qed

\vspace{5mm} \noindent Let us now assume $T$ and $M_{0}$ to be
countable. We introduce some more notation: Fix $k$, and let
$(Q_{k}^{i}:i<\omega)$ be an enumeration of the formulas in $Q_{k}$.
Let $\psi_{k}^{i}(x_{0},..,x_{k},y)$ be the formula
``$Q_{k}^{i}(x_{0},..,x_{k}) \wedge
\wedge_{j<k}(\phi(x_{j},y)\leftrightarrow \neg\phi(x_{j+1},y))$".
Let $\chi_{k}^{j,i}(y)$ be: $$``\exists
x_{0},..,x_{k}(\psi_{k}^{j}(x_{0},..,x_{k},y)\wedge (\neg\exists
x_{k+1}(\psi_{k+1}^{i}(x_{0},..,x_{k+1})))\wedge
\phi(x_{k+1},y))''.$$
\begin{Corollary}\label{3.9} For any $c\in {\bar M}$, $\phi(x,c)\in p$
if and only if there is $k \leq N_{\phi}$ and there is $i<\omega$
such that $c$ satisfies the formula $\chi_{k}^{j,i}(y)$ for all
$j<\omega$.
\end{Corollary}
\pf By Lemma \ref{3.8} and the notation.\qed

\vspace{5mm} \noindent Note that Corollary \ref{3.9} gives us an 
$F_{\sigma}$-definition $p$ over $M_{0}$.

In any case Theorem \ref{3.5} follows from Corollary \ref{3.9} as in
the proofs of Remark \ref{3.6}. Note that the only real assumption
on $p$ we need is that it is finitely satisfiable in some small
model (not necessarily $M_{0}$).\qed

\section{Groups with finitely satisfiable generics}
Here we introduce a certain desirable property of definable groups
which we call $fsg$ (standing for ``finitely satisfiable generics")
In Section 7 of the paper we  prove that definably compact groups
definable in $o$-minimal expansions of real closed fields have
$fsg$.

Again we fix a saturated model $\bar M$ of $T$. $G$ will denote a
group, definable in $ \bar M$ over $\emptyset$.

\begin{Definition}{\em  $G$ has $fsg$ ({\em finitely satisfiable
generics}) if there is some global type $p(x)$ and some small model
$M_{0}$ such that $p(x)\models ``x\in G"$, and every left translate
$gp=\{\phi(x):\phi(g^{-1}x)\in p\}$ of $p$ with $g\in G$, is
finitely satisfiable in $M_{0}$.}
\end{Definition}

The basic example of such a group is a {\em stable} group. (If $G$
is stable, then there exists a global generic type $p$ of $G$ in the
sense of stable group theory, namely every translate of $p$ does not
fork over $\emptyset$. But then by the characterization of forking
in the stable context, every translate of $p$ is finitely
satisfiable in any submodel $M_{0}$.) In simple theories, however,
definable groups will not, as a rule, have $fsg$. Also, the ordered
group $\la \mathbb{R}, <, +\ra$ does not have fsg. On the other hand
the {\em generically metastable} groups from $\cite{Hrushovski}$
which were introduced in connection with definability in
algebraically closed valued fields, {\em do} have $fsg$.

\vspace{5mm} \noindent For the remainder of this paper we call a
definable subset $X$ of $G$ (or the formula defining it) {\em left
generic} if finitely many left translates of $X$ cover $G$. Likewise
with right generic. $X$ is {\em generic} if it is both left and
right generic. A partial type $\Sigma(x)$ implying $x\in G$ is left
(right) generic if every formula in $\Sigma(x)$ is. Although this is
accordance with established vocabulary in the case of stable
theories, one should be aware that there is a discrepancy in the
case of simple theories. Notice that if $p$ is a global type in $G$,
and $X$ is a definable left generic subset of $G$ then some left
translate of $X$ (i.e. of the formula "$x\in X$") is in $p$.

\begin{Proposition} \label{4.2} Suppose that $G$ has $fsg$, witnessed 
by $p$ and
$M_{0}$, and let
$X\subseteq G$ be definable. Then
\newline
(i) $X$ is left generic iff $X$ is right generic (so we just say
generic).
\newline
(ii) $X$ is generic if and only if every left (right) translate of $X$
meets $M_{0}$.
\newline
(iii) $p$ is a generic type, as is any left or right translate of $p$.
\newline
(iv) If $X$ is generic and $X= X_{1}\cup X_{2}$ where the $X_{i}$ are
definable, then one of $X_{1}$, $X_{2}$ is generic.
\end{Proposition}
\pf Before we start let us note that
\newline
(*) $p^{-1}=\{\phi(x):\phi(x^{-1})\in p\}$ has the property that
every right translate of it is finitely satisfiable in $M_{0}$
\newline
(i) Suppose $X$ to be left generic. Then for any $c\in G$, $cX$ is
also left generic, so some left translate of $cX$ is contained in
$p$ whereby $cX $ is contained in some left translate $gp$ of $p$.
By the assumption (on $p$, $M_{0}$), $cX$ meets $M_{0}$, namely
there is $b\in G(M_{0})$ such that $b\in cX$, so $c^{-1}\in
Xb^{-1}$. We have shown that every element of $G$ lies in $Xb$ for
some $b\in G(M_{0})$. Compactness implies that finitely many right
translates of $X$ cover $G$, namely $X$ is right generic. The other
direction (right generic implies left generic) follows from (*) and
symmetry.
\newline
(ii) follows from the proof of (i).
\newline
(iii) If $X$ is in $p$ then every left translate of $X$ is in a left
translate of $p$ so meets $M_{0}$, whereby $X$ is generic by (ii).
\newline
(iv) If $X$ is generic, then $X$ is in a translate of $p$. Thus one
of $X_{1}, X_{2}$ is in the same translate of $p$. By (iii) one of
$X_{1}, X_{2}$ is generic.\qed

\vspace{5mm} Notice that Proposition \ref{4.2} implies that $G$ has
fsg, witnessed by $M_0$, if and only if every definable generic
subset of $G$ meets $M_0$ and the complement of every nongeneric set
is generic (the latter implies the existence of a generic type,
while the first implies that a generic type is finitely
satisfiable).

\vspace{5mm} \noindent It follows from (iv) that, assuming that $G$
has fsg, the set of nongeneric definable subsets of $G$ forms an
ideal $\cal I$ in the Boolean algebra of all definable subsets of
$G$. So for a definable subset $X$ of $G$, the stabilizer of $X$
modulo this ideal, namely $Stab_{\cal I}(X)=\{g\in G: gX\Delta X$ is
nongeneric\} forms a subgroup of $G$.  Note also that $Stab_{\cal
I}(X)$ is type-definable (by countably many formulas). On the other
hand for any global type $q$ of $G$, $Stab(q)$ is defined to be the
set of $g\in G$ such that $gq = q$. This is clearly a subgroup of
$G$ but on the face of it, has no definability properties.

\begin{Corollary}\label{4.3} Suppose that $G$ has $fsg$. Then
\newline
(i) There is a bounded number of (global) generic types,
\newline
(ii) $G^{00}$ exists.
\newline
(iii) For each (global) generic type $p(x)$, $Stab(p) = G^{00} =
\cap\{Stab_{\cal I}(X):X\in p\}$.
\end{Corollary}
\pf (i) Each generic type is finitely satisfiable in $M_{0}$ by
\ref{4.2} (ii). So there are a bounded number of them. (Any global 
type $p$
which is finitely satisfiable in a model $M_{0}$ is 
determined by $\{X\cap M_{0}:X\in p\}$.)
\newline
(ii) Let  (by part (i)) $\lambda$ be the number of global generic 
types of $G$. Fix a generic type $p$. Let $H$ be a type-definable 
subgroup of
$G$ of bounded index. So each coset of $H$ is in a translate of $p$.
The index of $H$ in $G$ is thus bounded by the number of (left)
translates of $p$, which is at most $\lambda$. So we have an absolute 
bound (independent of the monster model) on the index of 
type-definable subgroups of $G$ of bounded index, which clearly 
implies that $G^{00}$ exists.
\newline
(iii) Fix a global generic type $p$ of $G$. As $G^{00}$ has bounded index
some translate of $G^{00}$ is in $p$ (namely for some translate $C$ 
of $G^{00}$ $p(x)$ implies $x\in C$), whereby
\newline
(a) $Stab(p)\subseteq G^{00}$.
\newline
On the other hand clearly
\newline
(b) $\cap\{Stab_{\cal I}(X):X\in p\}\subseteq Stab(p)$, as $p$ only
contains generic definable sets.
\newline
So to conclude the proof of (iii) it suffices to prove
\newline
(c) For each definable $X\in  p$, $Stab_{\cal I}(X)\supseteq G^{00}$.
\newline
Suppose $X$ is defined over a small model $M$ containing $M_{0}$.
Note that if $g,h\in G$ and $tp(g/M) = tp(h/M)$ then $gX\cap G(M) =
hX\cap G(M)$, whereby $gX\Delta hX$ is not satisfiable in $M_{0}$
hence is nongeneric. It follows that the index of $Stab_{\cal I}(X)$
in $G$ is bounded by the number of types over $M$, that is to say,
$Stab_{\cal I}(X)$ has bounded index in $G$ hence contains $G^{00}$.
This proves (c) and completes the proof of the Corollary.\qed

\begin{Remark} \label{4.4} If $G$ has $fsg$ and $M_{0}'$ is any model then all
generic definable sets meet $G(M_{0}')$.
\end{Remark}
\pf Fix a formula $\phi(x,y)$ and $k<\omega$. By compactness there
is a finite subset $D$ of $G(M_{0})$ such that if $X\subseteq G$ is
defined by an instance of $\phi(x,y)$ and $k$ left translates of $X$
cover $G$, then $X$ meets $D$. Let $d$ be a finite tuple enumerating
$D$. Then the above property of $d$ can be expressed by a formula
without parameters. As this formula is realized in every model
$M_{0}'$ we are done.\qed

\vspace{5mm} \noindent The following will be helpful in carrying out
inductive proofs:
\begin{Proposition}\label{4.5} Let $G$ be a $\emptyset$-definable
group, and $N$ a $\emptyset$-definable normal subgroup of $G$.
Suppose that  $G/N$ and $N$ both have $fsg$. Then so does $G$.
\end{Proposition}
\pf Fix a small model $M_{0}$ witnessing that each of $G/N$ and $N$
have $fsg$. (In fact by 4.4 $M_{0}$ can be any submodel of ${\bar
M}$.)

We will freely use Proposition \ref{4.2}, applied to each of $G/N$
and $N$.

For $X$ a definable subset of $G$, let us define $Y_{X}$ to be
$\{g/N\in G/N: g^{-1}X\cap N$ is generic in $N$\}. Now $Y_{X}$ 
is
{\em not} necessarily definable, so we cannot apply directly 4.2. 
But $Y_{X}
= \cup_{i=1}^{\infty}Y_{X}^{i}$, where $Y_{X}^{i}$ is the set of
$g/N\in G/N$ such that $i$ left translates of $g^{-1}X\cap N$ by
elements of $N$ cover $N$. Each $Y_{X}^{i}$ {\em is} of course definable.

By compactness and the fact that $G/N$ is $fsg$ we have:
\newline
{\em Claim 1.} Finitely many left translates of $Y_{X}$ cover $G/N$
iff finitely many right translates of $Y_{X}$ cover $G/N$ iff  for
some $i<\omega$, $Y_{X}^{i}$ is generic in $G/N$.

\vspace{2mm}
\noindent
We will simply say ``$Y_{X}$ is generic in $G/N$" if the
equivalent conditions of Claim 1 hold.
\newline
{\em Claim 2.} If
$Y_{X}$ is generic in $G/N$ and $h\in G$ then each of
$Y_{hX}$ and
$Y_{Xh}$ are generic in $G/N$.
\newline
\pf Just notice that $Y_{X} = (h/N)^{-1}Y_{hX} = Y_{Xh}(h/N)^{-1}$.

\vspace{2mm}
\noindent
{\em Claim 3.} Suppose $X = X_{1}\cup X_{2}$ where the $X_{i}$ are 
definable. Then
\newline
(i) $Y_{X} = Y_{X_{1}}\cup Y_{X_{2}}$, 
and
\newline
(ii) If $Y_{X}$ is generic in $G/N$ then one of 
$Y_{X_{1}}$, $Y_{X_{2}}$ is generic in $G/N$.
\newline
\pf  (i) As $N$ has $fsg$, for each $g\in G$ $g^{-1}X\cap N$ is 
generic in $N$ if and only if
$g^{-1}X_{1}\cap N$ or $g^{-1}X_{2}\cap 
N$ is generic in $N$.
\newline
(ii) Assume that $Y_{X}$ is generic in 
$G/N$, so there are $h_{1},..,h_{n}\in G/N$ 
such that 
$\cup_{j=1,..,n} h_{j}Y_{X} = G/N$. By part (i) $G/N$ is covered by 
the $h_{j}Y_{X_{1}}$ 
together with the $h_{j}Y_{X_{2}}$ for 
$j=1,..,n$. Writing $Y_{X_{1}}$ as $\cup_{i<\omega}Y_{X_{1}}^{i}$ and 
likewise for $Y_{X_{2}}$ and applying compactness we see (as $G/N$ 
has $fsg$) that 
either some $h_{j}Y_{X_{1}}^{i}$ is generic in $G/N$ 
or some $h_{j}Y_{X_{2}}^{i}$ is generic in $G/N$. This suffices.

\vspace{2mm}
\noindent
{\em Claim 4.} If $Y_{X}$ is generic in $G/N$ then $X\cap
M_{0} \neq \emptyset$.
\newline
{\em Proof.} By Claim 1, let $i$ be such that $Y_{X}^{i}$ is generic
in $G/N$. Hence $Y_{X}^{i}\cap M_{0} \neq \emptyset$. This means
precisely that there is $h\in G(M_{0})$ such that $h/N\in
Y_{X}^{i}$. So $h^{-1}X\cap N$ is generic in $N$ and $h\in
G(M_{0})$. Now, since $N$ has fsg, the set $h^{-1}X\cap N$ contains
an element of $G(M_0)$, which clearly implies that $X$ does.

\vspace{2mm}
\noindent
To conclude the proof, for $X$ a definable subset of $G$, 
let us call $X$ *-generic if $Y_{X}$ is generic in $G/N$. By Claims 2 
and 3, the family of *-generics is closed under (left or right) 
translation, and the family of non *-generics forms a proper ideal. 
Hence there is a global *-generic type $p$ of $G$, and moreover by 
Claim 4, every translate of $p$ is finitely satisfiable in $M_{0}$. 
This shows that $G$ has $fsg$.  \qed

\begin{Remark}  The {\em fsg} property can also be formulated in terms of
measures. Namely, we say that a group $G$ {\em  has $fsg_m$} if
there is a Keisler measure $\mu$ on $G$ and there is some small
model $M_0$ of $\bar M$ such that for every $g \in G$, the measure
$g\mu$ (defined as $g\mu(X)=\mu(gX)$) is finitely satisfiable in
$M_0$.  It turns out that these formulations are equivalent:

If $G$ has {\em fsg} then the generic type gives the desired 0-1
measure. For the converse, one goes through the proof of \ref{4.2},
using the $fsg_m$ assumption instead of {\em fsg}.
\end{Remark}

One could ask if the results of this section hold under the weaker assumption
that there exists a type over a large saturated model, with a small number of
left translates.  However variants of the following example will show 
that various lemmas
can fail, even for rank one simple theories with a translation invariant type.

\begin{Example}
Let $T_0$ be the theory of vector spaces over $GF(2)$ with a symmetric
irreflexive binary relation R (bearing no particular relation to $+$.)
Let $T$ be the model completion of $T_0$, $M \models T$.  Let
  $C \subseteq M$ and $D \subseteq M \setminus \{0\}$ be
  arbitrary subsets.
Then $\{R(x+b,c):  c \in C\} \cup\{ \neg R(x+b,c):  c \notin C\} \cup 
\{ R(x+c,x+b): b-c \in D\} \cup \{ \neg R(x+c,x+b): b-c \notin D\}$ 
determines a complete type, which is $M$-translation invariant.
But there are essentially no generic formulas.
Theorem \ref{3.5}.

\end{Example}

\section{Definably amenable groups}

It is a convenient time to introduce the notion ``definable
amenability". Recall that an abstract (or discrete) group $G$ is
said to be {\em amenable} if there exists a left invariant finitely
additive probability measure on the family of {\em all} subsets of
$G$. Any solvable group is amenable.

\begin{Definition} Let $G$ be a definable group. We call $G$
definably amenable if there is a left invariant Keisler measure on
$G$.
\end{Definition}

\begin{Remark} (i) Any amenable group is definably amenable.
\newline
(ii) Suppose $T$ has a model $M_{0}$ such that $G$ is defined over
$M_{0}$ and $G(M_{0})$ has a compact (Hausdorff) group topology such
that every definable subset of $G$ is Haar measurable. Then $G$ is
definably amenable.
\newline
(iii) If $K$ is a (saturated) algebraically closed valued field and
$n>1$ then $SL(n,K)$ is not definably amenable.  
\newline (iv) If 
$R$ is an expansion of a real closed field,
then $PSL(2,R)$ is not definably amenable.   
\newline
(v) 
$SO(3,\Rr)$ is definably amenable,
but not amenable as a pure group.
\end{Remark}
\pf (ii) is proved by construction (*) applied to the (unique)
normalized Haar measure on $G(M_{0})$.
\newline
(iii) This follows by a similar proof to that in \cite{Hrushovski}
showing that $SL(n,K)$ has no definable left generic type.
(iv) Suppose $\mu_1$ is a left invariant Keisler measure  on $PSL(2,R)$.
Recall the transitive action of $SL(2,R)$ on $\Pp^1(R) = R \union \{\infty\}$.
Define a Keisler measure $\mu$ on $\Pp^1(R)$ by
$\mu(X)  = \mu_1(\{g \in PSL(2,R): g \cdot 0 \in X \}$.  Then $\mu(hX) =
\mu_1(\{g \in PSL(2,R): h^{-1}g \cdot 0 \in X \} = \mu(X)$.  But let 
$U$ be a small ball
around $0$; then using inversion we find $gU$, a   ball around $\infty$; while
using multiplication we can find $hU$ such that $\Pp^1(R) = gU \union hU$.
So $\mu(U) \geq 1/2 \mu(\Pp^1(R))$.  This is true for an arbitrary 
small ball around any point,
e.g. $0,1,\infty$, giving $3/2 \leq 1$,
a contradiction.  (v)  Every definable set is Lebesgue measurable. 
The pure group  statement
is due to Hausdorff, Banach and Tarski, see below.

\vspace{5mm} \noindent Before continuing we take the opportunity to
give a characterization of definable amenability (and the
construction of an invariant Keisler measure on $G$ from a suitable
Grothendieck group of $G$). Fix a definable group $G$. By a
nonnegative cycle in $G$ we mean a ``finite disjoint union" of
definable subsets of $G$. Notationally consider a nonnegative cycle
as $\{k_{1}X_{1},..,k_{n}X_{n}\}$ where the $k_{i}$ are nonnegative
integers, and $X_{i}$ are pairwise distinct definable subsets of
$G$. There is the obvious notion of a map between two such cycles
being definable, $1-1$ and given by piecewise left translations.
\begin{Definition} A definable paradoxical decomposition of $G$ is a
definable $1-1$ piecewise translation from the disjoint union of $G$
and $Y$ to $Y$, for some nonnegative cycle $Y$.
\end{Definition}

This is a variant of a notion due to Hausdorff.  He actually 
considered a stronger notion,
that would rule out the existence of any invariant finitely additive 
measure, not necessarily
non-negative.  His construction of a paradoxical decomposition of the 
two-dimensional sphere, or of $SO(3,\Rr)$, completed by Banach and 
Tarski, requires the axiom of choice and
is not represented in a definable way.  On the other hand we have:

\begin{Proposition} \label{amenability} $G$ is definably amenable if 
and only if $G$ does
not admit a definable paradoxical decomposition.
\end{Proposition}

\vspace{3mm} \noindent Before entering the proof we introduce the
relevant Grothendieck (semi-) group. Let $K_{semi}(G)$ be the
semigroup whose elements are the nonnegative cycles
$\sum_{i}k_{i}X_{i}$ in $G$
   modulo the equivalence relation of being in definable bijection by
piecewise left translations. A typical element of $K_{semi}(G)$ can
be written in the form $k_{i}[X]_{semi}$ where $[X]_{semi}$ is the
class of the definable set $X$ in $K_{semi}(G)$. Addition in the
semigroup is the obvious thing.

Let us make a further identification: let $x_{1}, x_{2}\in
K_{semi}(G)$. Define $x_{1}\sim_{0} x_{2}$ if there is $y\in
K_{semi}(G)$ such that $x_{1} + y = x_{2} + y$. Then the collection
of $\sim_{0}$-classes, together with formal inverses, constitutes
the Grothendieck group $K_{0}(G)$. The class in $K_{0}(G)$ of a
definable subset $X$ of $G$ is denoted $[X]_{0}$. (Likewise for a
nonnegative cycle $Y$.)

\vspace{5mm} \noindent {\em Proof of Proposition 5.4.} Suppose first
that $\mu$ is a left invariant Keisler measure on $G$. Then a
definable paradoxical decomposition could not exist since then we
would have $\mu(G) + \mu(Y) \leq \mu(Y)$, contradicting $\mu(G) =
1$.
\newline
Conversely suppose $G$ has no definable paradoxical decomposition.
Let $P_{0}$ be the subsemigroup of $K_{0}(G)$ generated by the sets
$[X]_{0}$ where $X$ is definable.
\newline
{\em Claim.} $-n[G]_{0}\notin P_{0}$ for all $n> 0$.
\newline
{\em Proof.} Otherwise $-n[G]_{0} = [Y]_{0}$ (in $K_{0}(G)$) for
some nonnegative cycle $Y$. But then $n[G]_{0} + [Y]_{0} = 0$ in
$K_{0}(G)$, so $n[G]_{semi} + Y_{semi} + [Z]_{semi} = [Z]_{semi}$ in
$K_{semi}(G)$ for some nonnegative cycle $Z$. But then clearly there
is a definable injective piecewise translation map from the disjoint
union of $G$ and $Z$ into $Z$, contradicting our assumption.

Let $B$ be the tensor product of $\Q$ with $K_{0}(G)$, and $$P =
\{\alpha x: \alpha\in {\Q}, \alpha > 0, x\in P_{0}\}$$ By the Claim
$-[G]_{0}\notin P$.
Let $P'$ be a
maximal subset of $B$  containing $P$, closed under multiplication by 
positive rationals
and addition, and
  such that $-[G]_{0}\notin P'$.  Define  a partial ordering on $B$: 
$x \leq y \iff y-x \in P'$.

  {\em Claim.} $\leq$ is a total ordering on $B$.
  \newline
{\em Proof.} We have to show that for any $a \in B$, either $a \in 
P'$ or $-a \in P'$
If $a \notin P'$, let $P'' = \{x+\alpha a:  x \in P', \alpha \in \Qq, 
\alpha>0 \}$.
Then by maximality $-[G]_0 \in P''$, i.e. $-[G]_0 =x + \alpha a, x \in P$,
so $-a = \alpha ^{-1} ([G]_0 +   x) \in P'$.

Now exists a  unique order preserving semigroup
homomorphism $h: B \to \R^{\geq 0}$ such that $[G]_{0}$ goes to $1$.
(Namely $h(b) = \alpha$ iff
$[b:[G]_0] = [\alpha:1]$ in the sense of \cite{euclid} V. Def. 5, i.e.
for all $m,n \in \Nn$, $mb < n  [G_0]$ iff $m \alpha < n \cdot 1$).
Let $\mu(X) = h([X]_0)$; this is clearly a  left invariant Keisler 
measure on $G$.

  \paragraph{Definable amenability of volumes.}

We will see in \S 8 that definably compact groups in o-minimal 
structures are definably
amenable.  The proof uses a deep structure theory for such groups. 
In the following paragraphs, not otherwise used in this paper, we 
consider a similar amenability property of
definable compact definable sets.   We do not know if this property 
holds in all o-minimal theories, but when it does we give a soft 
proof of definable amenability of definable compact groups.  In 
particular this is valid for  o-minimal
expansions of RCF that are finitely satisfiable in expansions of 
$(\Rr,+,\cdot)$.

The proof actually yields more:  that any definable group $G$, not 
necessarily definably compact, is {\em definably amenable for compact 
sets}.    By definition this means:
let $Def_{bdd}(G)$ be the family of definable subsets of definably 
compact subsets of $G$.
Then
for any $X \in Def_{bdd}(G)$  with nonempty interior,  there exists a 
translation invariant finitely additive
   $\mu: Def_{bdd}(G) \to \Rr^{\geq 0} \union \{\infty\}$ with $\mu(X)=1$.


Let $T$ be an o-minimal expansion of RCF, and fix $n \geq 1$.  By 
``almost all' we will mean:
away from a definable set of dimension $<n$.  If $f: R^n \to R^n$ is 
definable, $|J f|(c)$ denotes
the absolute value of the determinant of the matrix of partial 
derivatives of $\phi$ at $c$;
it exists almost everywhere.
Let $V[n]$ be the set of bounded definable functions $R^n \to R^{\geq 
0}$  with bounded support.  By an {\em isomorphism} $\phi: f \to g$
we mean a definable bijection $\phi$ from a definable set containing 
the support of $f$
to one containing the support of $g$, such that $f(x)= |J \phi |(x) 
\cdot g(\phi(x)) $
almost everywhere. More generally define $f \sim g$ if one can write
$f=\sum_{i=1}^n f_i, g=\sum_{i=1}^n g_i$ with $f_i,g_i$ isomorphic. 
Let $[f]$ denote the
$\sim$-class of   $f$, and let
$K_{semi}(V[n]) = \{[f]: f \in V[n]\}$.  Define $[f]+[g] = [f+g]$. 
Let $K(V[n])$ be the
corresponding group.  We say that $T$ is {\em definably amenable for 
volumes} if for each $n$
and any $f \in V[n]$, either $f=0$   a.e., or there exists an 
order-preserving semigroup homomorphism $I: K_{semi}(V[n]) \to 
\Rr^{\geq 0} \union \{\infty\}$ with $0 < I(f) < \infty$.

\begin{Proposition}

Conditions (1) , (2) are equivalent.  Also (3) or (4) imply (1) and (2);
while (1) or (2)  implies (5) and (5) implies (6).

   \begin{enumerate}

   \item $T$ is amenable for volumes.
   \item If $f,h \in V[n]$ and $f \sim f+h$, then $h=0$ a.e.

     \item Every finite $T_0 \subseteq T$ has a complete archimedean model.
    \item  $T$ has definable primitives:  for
every definable function $f: R \to R$ in a model of $T$, there exists 
a definable function $F$
    such that almost everywhere $F'=f$.

    \item Every definable group $G$ of $T$ is definably amenable for 
compact sets.
     \item Every definably compact group in $T$ is definably amenable.
\end{enumerate}
\end{Proposition}

\proof The proof of equivalence of (1) and (2) is identical to the 
proof of Proposition
\ref{amenability}.

If (2) fails, then there exists a finite $T_0$ describing the 
situation.  For instance
if there exists an isomorphism $\phi: f \to f+h$, then $\phi$ is differentiable
away from a set $Y$ of dimension $n-1$; there exist $W_1,\ldots,W_k 
\subseteq R^{n-1}$
and continuously differentiable $e_i: W_i \to Y$ with $\union_{i=1}^k 
e_i(W_i) = Y$;
and $f+h= |J(\phi)|  \cdot f \circ \phi $ away from $Y$.  $T_0$ can 
state this, as well
as the   boundedness and piecewise differentiability  of $f,h$, and 
that $h>0$ on some open set.
   If $T_0$ has a complete archimedean model $N=(\Rr,+,\cdot, 
f^N,h^N,\ldots)$, then $f^N,h^N$ are bounded  integrable functions, 
and $\phi^N$ is $C^1$ away from $Y^N$, so by the
   change of variable formula $\int f^N = \int (f+h)^N$ (where $\int$ 
is the Lebesgue or Riemann integral.)   But also $\int h^N > 0$, a 
contradiction.    This shows that (3) implies (2).  The proof that(4) 
implies (2) is similar:  by a compactness argument, a function 
$f(x_1,\ldots,x_n)$
   has a definable primitive $F_1(x_1,\ldots,x_n)$ with respect to the 
first variable,
   i.e. $\partial F_1 / \partial x_1 = f$ a.e.  Now one can define 
integration using iterated
   integrals, and prove the change of variable formula and additivity 
using o-minimality.
The proof of (3) implies (2) used no more than this.

(6) is obviously a special case of (5).

To prove (5) from (1) let $n=\dim(G)$.
  Fix an
identification of some neighborhood of $1$ in $G$ with an open 
neighborhood of $0$ in $R^n$.
Let $K_0$ be the set of subsets $Y$ of $G$
contained in $b(int(U))$ for some   injective continuously 
differentiable definable map  $b: U \to G$,
$U$  a definably compact subset of $R^n$ with interior $int(U)$. 
We begin by defining a map $\psi: K_0 \to K_{semi}(V[n])$.

For $g \in G$,  let $T_g: G \to G$, $T_g(x) = g^{-1}x$.
Given $Y \in K_{0}$, find a definably compact set $U \subseteq R^n$ 
and a definable injective
  $C^1$ map $b: U \to G$, with $Y \subseteq b(int(U))$.
   Let $f(x)=0$ for $x \notin b ^{-1} (Y)$, and for $x \in b ^{-1}(Y)$
let $f(x) = |Jg|(x)$, where $g(y) = T_{b(x)} \circ b$.
(Here we use the identification of a neighborhood of $1$ with a 
neighborhood of $0$ in $\Rr^n$;
so $g: U \to \Rr^n$, and the Jacobian $Jg$  is defined.)
By continuity and definable compactness, $f$ is bounded on $U$.
If we pick a different $b': U' \to Y$, with corresponding $g',f'$,
then $b' = b \circ e$ for some $e: U' \to U$ (defined on a 
neighborhood of the support of $b$),
namely  $e(u')=b ^{-1}(b'(u'))$  on $(b')^{-1}( (b(int(U))$.  We have
$g' = g \circ e$, $|Jg'|(x) = |Jg|(e(x)) |Je|(x)$
so that $f$ is isomorphic to $f$ and $[f]=[f'] \in K_{semi}(V[n])$. 
Hence $[f]$ does
not depend on the choice of $(U,b)$ and we can
define $\psi(Y) = [f]$.

Given $h \in G$, let $b'' = hb$; then
$T_{b''(x)} \circ b'' = T_{b(x)} \circ b$, so $\psi(hU) = \psi(U)$. 
Thus $\psi$ induces a
well defined map $K_1 \to K_{semi}(V[n])$, where
  $K_1 = \{[Y]: Y \in K_0 \} \subseteq K_{semi}(G)$.  It is clear that 
$\psi(a+b)=\psi(a)+\psi(b)$
  when $a,b,a+b \in K_1$, and that $a+b \in K_1 $ implies $a \in K_1$. 
It follows
  that $\psi$ extends to
homomorphism of ordered   semigroups $\sum K_1 \to K_{semi}(V[n])$,
where $\sum K_1$ is the semigroup generated by $K_1$.

According to \cite{BO2}, for any definably compact $Z \subseteq G$ 
there is a $C^1$ group manifold structure on $G$ with finite
chart $\{b_i: W_i \to G: i=1,\ldots,r\}$ (with the $W_i$ open subsets 
of $R^n$) and closed bounded
  $U_i \subseteq W_i$, such that $X \subseteq \union_{i=1}^r 
b_i(int(U_i))$.  Hence $\sum K_1 = K_{semi}(G)$.

   If $X$ has nonempty interior, then $\psi(X)$ cannot vanish almost everywhere,
so by (1) there exists homomorphism   $\mu: K_{semi}(V[n]) \to 
\Rr^{\geq 0} \union \{\infty\}$
with $\mu(\psi(X))=a >0$.   Now
   $(1/a) \mu \circ \psi$ demonstrates (5).

It is of course possible to combine (3) and (4), i.e. it suffices 
that every finite $T_0 \subseteq T$
be extendible to an $o$-minimal theory with definable primitives, or 
to one with an archimedean  model.

\begin{Question}  Is every o-minimal theory amenable for volumes? 
\end{Question}

\section{Groups with NIP}
Here we concentrate on definable groups in theories with NIP.

Suppose  that $\mu$ is a Keisler measure on a definable group $G$.
Then for any $g\in G$, we have another Keisler measure $g\mu$ on
$G$, namely $g\mu(X) = \mu(gX)$. We say that $\mu$ is {\em left
invariant} if $g\mu = \mu$ for all $g\in G$. Likewise for right
invariant. The existence of a left invariant {\em type} of $G$ is a
very strong property. For example if $G$ is stable, this implies
that $G$ is connected and the left invariant type is the unique
generic type of $G$. However, even if there is NO invariant type,
one may hope for there to exist an invariant measure.

\vspace{5mm} \noindent The next proposition, due to Shelah
\cite{Shelah}, gives the existence of $G^{00}$ for any definable (or
even type-definable) group $G$ in a theory with NIP. We had
originally proved this under the additional assumption that $G$ was
definably amenable. In any case thanks to Shelah for allowing us to
include the result and a proof.
\begin{Proposition}\label{5.5} Assume $T$ has NIP. Let $G$ be a 
definable group in
$\bar M$, defined over $\emptyset$ say. Then $G$ has a smallest
type-definable subgroup of bounded index. If $G^{00}$ is such then it is
type-definable over $\emptyset$ and has index at most $2^{|T|}$.
\end{Proposition}
\pf It is easy to see that any type-definable subgroup of $G$ is
the intersection of a family of subgroups each of which is
type-defined by countably many formulas (see for example Remark
1.4(ii) in \cite{BOPP}). So it suffices to prove that any subgroup
$H$ of $G$ which is type-defined by countably many formulas and has
bounded index in $G$ has only a bounded number of distinct
conjugates (under automorphisms of the ambient structure). So let us
suppose, for a contradiction, that $a$ is a countable tuple, $H_{a}$
is type-definable by a countable partial type $\Sigma(x,a)$ over
$a$, $H_{a}$ is a subgroup of bounded index in $G$, and that
$\{H_{a'}: tp(a') = tp(a)\}$ is unbounded (where $H_{a'}$ is
type-defined by $\Sigma(x,a')$). So by Erdos-Rado we have some
indiscernible sequence $\la a_{i}:i<\omega\ra$ of realizations of $p
= tp(a)$ such that $H_{a_{i}} \neq H_{a_{j}}$ for $i\neq j$.
\vspace{2mm}

\noindent {\em Claim 1.} Fix $i_{0} <\omega$. Then $\cap\{H_{a_{j}}:
j<\omega, j\neq i_{0}\}$ is NOT contained in $H_{a_{i_{0}}}$.

\noindent {\em Proof of claim 1.} Suppose otherwise. We can
``stretch" the indiscernible sequence $\la a_{i}:i<\omega\ra$ by
inserting some $(b_{\alpha}:\alpha < \kappa)$ in place of
$a_{i_{0}}$ (for any $\kappa$). But then each $H_{b_{\alpha}}$
contains $\cap_{j\neq i_{0}}H_{j}$. But $\alpha\neq \beta$ implies
$H_{b_{\alpha}} \neq H_{b_{\beta}}$. So for any $\kappa$ we can find
at least $\kappa$ many distinct subgroups of $G$ each of which
containing $\cap_{j\neq i_{0}}H_{j}$. As the latter has bounded
index in $G$, we get a contradiction, proving the claim.

\vspace{2mm} \noindent The claim clearly applies also to any
stretching $\la a_{\alpha}\ra $ of the indiscernible
sequence $\la a_{i}:i<\omega\ra$. So for each $\alpha$, let
$c_{\alpha}$ be such that $c_{\alpha}\notin H_{\alpha}$ but
$c_{\alpha}\in H_{\beta}$ for all $\beta\neq \alpha$. Again by
Erdos-Rado we may assume that the sequence $\la
(a_{\alpha},c_{\alpha}):\alpha <\kappa\ra$ is indiscernible.
\newline
We may assume $\Sigma(x,a) =\{\phi_{n}(x,a):n<\omega\}$
where $n<m$
implies $\models\phi_{m}(x,a)\rightarrow \phi_{n}(x,a)$.

\vspace{2mm}
\noindent
{\em Claim 2.} There is $n<\omega$ such that for any $\alpha$ and any
$d_{1},d_{2}\in H_{\alpha}$,
\newline
$\models \neg\phi_{n}(d_{1}\cdot c_{\alpha}\cdot d_{2}, a_{\alpha})$.
\newline
{\em Proof of claim 2.} As $tp(a_{\alpha},c_{\alpha})$ does not depend on
$\alpha$, it is enough to prove it for a fixed $\alpha$.
As $c_{\alpha}\notin H_{\alpha}$ we have the implication:
\newline
$y_{1}, y_{2}\in H_{\alpha} \models \vee_{n}\neg\phi_{n}(y_{1}\cdot
c_{\alpha}\cdot y_{2},a_{\alpha})$. Now apply compactness.

\vspace{2mm} \noindent We may clearly assume $n=0$ in Claim 2.
\vspace{2mm}

\noindent{\em Claim 3} For each finite $w\subset \kappa$ there is
$d_{w}$ such that for all $\alpha$, $\models
\phi_{0}(d_{w},a_{\alpha})$ iff $\alpha\notin w$.

\noindent {\em Proof of claim 3.} Let $d_{w}$ be the product of the
$c_{\beta}$ for $\beta\in w$. So if $\alpha\notin w$, then as
$c_{\beta}\in H_{\alpha}$ for each $\beta\in w$, $d_{w}\in
H_{\alpha}$ hence satisfies $\phi_{0}(x,a_{\alpha})$. On the other
hand if $\alpha\in w$ then we can write $d_{w}$ as $d_{1}\cdot
c_{\alpha}\cdot d_{2}$ where $d_{1}, d_{2}\in H_{\alpha}$ (by an
argument as above). So then we apply Claim 2.

\vspace{2mm} \noindent Claim 3 shows that $T$ has the independence
property, a contradiction. So $G^{00}$ exists. Its type-definability
over $\emptyset$ follows by uniqueness (any type-definable set which 
is $\emptyset$-invariant is type-definable over $\emptyset$, by 
quantifying out the parameters and using saturation). The bound on 
the index is
clear too.\qed

\vspace{5mm} \noindent The existence of $G^{00}$ (for $G$ a
definable group in a saturated model of $T$) had been proved earlier
in various special cases. For example for $o$-minimal theories in
\cite{BOPP}. In fact the latter proved in addition that $G/G^{00}$
is a compact Lie group. For groups definable over ${\Q}_{p}$ in a
model of $Th(({\Q}_{p})_{an})$ this was done in
\cite{Onshuus-Pillay}. For groups definable in Pressburger
Arithmetics, it follows from work of Onshuus \cite{Onshuus}.

\vspace{5mm} \noindent Here is an  application of Proposition
\ref{5.5}. Let us fix a compact Hausdorff group $\la G,\cdot,...\ra$
equipped with additional first order structure. We use the term $G$
to also denote this structure. Let us assume that (i) $Th(G)$ has
the NIP, (ii) any definable subset of $G$ is Haar measurable (with
respect to the unique normalized Haar measure on $G$), and (iii)
there is a neighbourhood basis of the identity of $G$ consisting of
definable sets, say $U_{i}$ for $i\in I$.

Let $G^{*}$ be a saturated elementary extension of $G$. So
$\cap_{i\in I} U_{i}^{*}$ is the group of ``infinitesimals", denoted
by $inf(G^{*})$ of $G^{*}$, and the quotient group (with the logic
topology) is precisely $G$. By Proposition \ref{5.5}, $(G^{*})^{00}$
exists, and in fact we have:
\newline
{\em CLAIM.} $(G^{*})^{00}$ is precisely the group
$inf(G^{*})$ of infinitesimals of $G^{**}$.
\newline
{\em Proof.} By \ref{5.5}, $(G^{*})^{00}$ is type-definable over
$\emptyset$. As we already know that $inf(G^{*})$ is type-definable,
and of bounded index, it suffices to prove that any subgroup $H$ of
$G^{*}$ which is
    type-definable over $G$ by a countable set of formulas,
and has bounded index, contains $inf(G^{*})$. Let $H$ be
such, and suppose
$H =
\cap_{n}X_{n}$, where $X_{n}$ is definable over $G$. We may assume that
$X_{n}^{-1}\cdot X_{n}\subseteq X_{n-1}$ for all $n>0$. Fix $n$. As $H$
has bounded index in
$G^{*}$, finitely many translates of $X_{n}(G)$ cover $G$ whereby the
Haar measure of $X_{n}(G)$ is $>0$.

It follows (cf. the  chapter on convolutions in \cite{HewittRoss}) 
that $(X_{n}^{-1}\cdot X_{n})(G)$ has
interior in $G$, so $(X_{n-1}^{-1}\cdot X_{n-1})(G)$ contains an
open neighbourhood of the identity of $G$. Thus $X_{n-2}(G)$
contains some $U_{j}$. Hence
   $H$ contains $inf(G^{*})$, and the claim is proved.

\vspace{5mm}
\noindent
Now measures come back into the picture.
The following was proved in the stable case in
\cite{Newelski-Petrykowski}.
\begin{Proposition} \label{5.6} Suppose $T$ has NIP, and $G$ is a
$\emptyset$-definable group in $\bar M$ with $fsg$. Then there is a
left invariant Keisler measure $\mu$ on $G$, which is moreover
finitely satisfiable in some small model $M_{0}$.
\end{Proposition}
\pf We will use \ref{4.2} and \ref{4.3}. Let us fix a global generic
type $p$ of $G$ over ${\bar M}$, such that $p(x)$ implies $x\in
G^{00}$. The measure we construct will depend on $p$. We will first
prove the proposition in the case where $T$ is countable. By Remark
\ref{4.4} let us fix a countable model $M_{0}$ such that all generic
definable subsets of $G$ meet $G(M_{0})$.  Let $\textbf{m}$ be the
(unique) normalized Haar measure on the compact group $G/G^{00}$.
\vspace{2mm}

  \noindent{\em Claim.} Let $X\subseteq G$ be definable.
Then
\newline
(i) For $g\in G$, whether or not $gX\in p$ depends only on the coset of
$g$ modulo $G^{00}$.
\newline
(ii) $\{g/G^{00}:gX\in p\}$ is a Borel subset of $G/G^{00}$ (so is Haar
measurable).
\newline
{\em Proof.} (i) follows because $Stab(p) = G^{00}$.
\newline
(ii): By Remark \ref{4.4} let $M_{0}$ be a countable model such that 
$X$ is over $M_{0}$ and
all generic definable subsets of $G$ meet 
$M_{0}$. In particular $p$ is finitely satisfiable in $M_{0}$.
By 
Corollary 3.9, there are partial types $\Psi_{i}(y)$ over $M_{0}$ for 
$i<\omega$
such that for $g\in G$, $gX\in p$ iff $\models 
\vee_{i<\omega}\Psi_{i}(g)$. Let $C_{i}$ be the closed subset of 
$G/G^{00}$ determined by $\Psi_{i}(y)$, namely the image of the 
solution set of $\Psi_{i}$ under the natural map taking to $G$ to 
$G/G^{00}$. Then by part (i) of the Claim, $\{g/G^{00}:gX\in p\}$ is 
precisely $\cup_{i}C_{i}$, hence Borel. 

\vspace{2mm} \noindent By the Claim, we can define $\mu_{p}(X) =
{\textbf{m}}(\{g/G^{00}:gX\in p\})$. Then $\mu_{p}$ is finitely
additive. For left invariance: let $g'\in G$, then $\{g/G^{00}:g\in
g'X\} = \{g/G^{00}:g\in X\}g'/G^{00}$, so by right invariance of
${\textbf{m}}$, $\mu_{p}(g'hX) = \mu_{p}(X)$.

Finally let us note that $\mu_{p}(X) > 0$ if and only if $X$ is generic.
Right implies left is  true by invariance of $\mu_{p}$. But if $X$ is
nongeneric, then no translate of $X$ is in $p$, so $\{g/G^{00}:gX\in p\}$
is empty, hence $\mu_{p}(X) = 0$.

As we already know that every generic definable subset of $G$ meets
$G(M_{0})$ for some small model $M_{0}$ we obtain finite 
satisfiability of $\mu_{p}$
in $M_{0}$.

\vspace{2mm} \noindent So we have proved the proposition when $T$ is
countable. For the general case: given a definable subset $X$ of
$G$, let $L_{0}$ be a countable sublanguage of $L$ in which $G$ and
$X$ are definable. Let $p_{0}$ be the reduct of $p$ to $L_{0}$ and
let $G_{0}^{00}$ be the smallest $L_{0}$-type-definable subgroup of
$G$ of bounded index. Let $f$ be the canonical surjective 
homomorphism from 
$G/G^{00}$ to $G/G_{0}^{00}$. Clearly $f$ is 
continuous. Let $U = \{g/G^{00}:gX\in p\}$ and $U_{0} = 
\{g/G_{0}^{00}: gX\in p_{0}\}$. Then $U = f^{-1}(U_{0})$. But by the 
Claim in the countable case, $U_{0}$ is Borel. Hence $U$ is also 
Borel. So the Claim holds in general, and as above we obtain our 
measure $\mu_{p}$. \qed

\vspace{2mm}
\noindent
It is natural to ask whether the measure $\mu_p$ defined above indeed
depends on the type $p$ or not. This and related issues will be 
tackled in a subsequent paper.

  \vspace{5mm} \noindent Our final result of this section will
provide in a sense the missing link in the proof of the $o$-minimal
conjectures.
\begin{Proposition}\label{5.7} Suppose that $T$ has NIP, and $G$ is a group
definable in ${\bar M}$ such that $G$ is definably amenable, and
the set ${\cal I}$ of non (left) generic definable subset of $G$ forms
an ideal. Then
\newline
(i) there are only a bounded number of definable subsets of $G$ modulo
the equivalence relation $X\sim_{\cal I} Y$ ($X\Delta Y\in {\cal I}$).
\newline
(ii) for each definable left generic
$X\subseteq G$,
$Stab_{\cal I}(X)$ $(=\{g\in G: gX\Delta X$ is nongeneric\}) is a
(type-definable) subgroup of bounded index.
\end{Proposition}
\pf Let $\mu$ be a left invariant Keisler measure on $G$. Note that
if $X$ is a left generic definable subset of $G$ then $\mu(X)
>0$ (as finitely many left translates of $X$ cover $G$ and these have all
the same $\mu$-measure as $X$). So if there unboundedly many
$\sim_{\cal I}$-classes there will also be unboundedly many
$\sim_{\mu}$-classes, contradicting Corollary \ref{3.4}. This proves
(i).
\newline
(ii) follows immediately.\qed

\section{Interlude: Ind-definable and locally compact groups}
As one of the authors remarked ``it seems a pity to lose $SL_{2}(\R)$".
So we give the notion of an Ind-definable group, point out that
quotienting by a type-definable normal subgroup of bounded index yields a
locally compact group, and develop analogues of some of the results so
far for Ind-definable groups.
We also state an appropriate version of the $o$-minimal conjectures from
\cite{Pillay}. In any case we will be brief.

We still work in a saturated model ${\bar M}$. Ind-definable stands
for ``inductive limit of definable sets". For notational reasons we
will take the index set to be $\N$. So an {\em Ind definable set}
$X$ will be by definition a sequence $(X_{n}:n\in \N)$ of definable
sets together with definable embeddings $f_{n}: X_{n}\to X_{n+1}$
for $n\in N$. The points of $X$ correspond to sequences
$(x,f_{n_0}(x), f_{n_0+1}(f_{n_0}(x)),...)$ for some $x\in
X_{n_{0}}$ and $n_{0}\in \N$. It is convenient to view the $f_{n}$
as inclusion maps, and so $X$ as the increasing union
$\bigcup_{n}X_{n}$. There are natural notions of an Ind-definable
relation on $X$ and Ind-definable functions between Ind-definable
sets. For example an Ind-definable function $g$ between $X =
\bigcup_{n}X_{n}$ and $Y = \bigcup_{n}Y_{n}$ is a function from $X$
to $Y$ such that for every $m,n$ $\{x\in X_{m}:g(x)\in Y_{n}\}$ is
definable and the restriction of $g$ to this set is definable. We
also have the obvious notion of an Ind-definable set, function,..
being defined over a given set $A$ of 
parameters.

\begin{Definition}{\em  An {\em Ind-definable group} $G$ 
is something
of the form $G = \bigcup_{n}G_{n}$ where
$G_{n}$ are definable sets, $m:G\times G\to G$ is a group operation
and when restricted to $G_{n}\times G_{n}$ has values in $G_{n+1}$
(and is definable), and inversion when restricted to $G_{n}$ has
values also in $G_{n}$.}
\end{Definition}

We could also say that an Ind-definable group $G$ is a group object 
in the category
of Ind-definable sets, noting that up to isomorphism 
$G$ has the explicit form given in Definition 7.1.

By a {\em 
definable} subset of $G$ we mean a definable subset of
some $G_{n}$. Likewise a complete type extending $G$ will be
``concentrate" on some $G_{n}$.

For various reasons we will assume that
\newline
(*) $G_{0}$ generates $G$ as a group.

\vspace{5mm} \noindent{\bf Examples} A basic example we have in mind
for an Ind-definable group is a subgroup of a definable group $G$
that is generated by a definable set $G_0\subseteq G$ (such groups
were called in \cite{Peterzil-Starchenko1}, ``$\bigvee$-definable
groups'' and in \cite{Edmundo} ``locally definable''). Another is
the universal cover of $\la [0,1),+(mod 1)\ra$, obtained as an
increasing union of intervals $[-,n,n]$ and the obvious group
operation. The group of definable automorphisms of a definable group
$G$, say in a countable language, can also be viewed as an
Ind-definable group, where the $G_n$'s in the definition are
obtained via the various definable families of automorphisms of $G$.
Finally, ``an infinite dimensional'' example is, for a definable
group $G$, the increasing union of $G, G\times G,\ldots,
G^n,\ldots$, with the group operation acting coordinate-wise (such
spaces are called by A. Piekosz, in preliminary notes, ``weakly
definable spaces'').

Here are some analogues of the basic notions:
\begin{Definition}{\em  Let $G = \cup_{n}G_{n}$ be an Ind-definable group.
\newline
(i) Let $X$ be a definable subset of $G$ (i.e. of some $G_{n}$). We
call $X$ {\em left generic in $G$} if for each $m$ finitely many
left translates of $X$ by elements of $G$ cover $G_{m}$.
\newline
(ii) By {\em a type-definable subgroup of $G$} we mean a subgroup
$H$ of $G$ which is at the same time a type-definable subset of some
$G_{n}$.
\newline
(iii) By {\em a Keisler measure on $G$} (or on any Ind-definable set
for that matter) we mean a finitely additive real-valued function
$\mu$ on definable subsets of $G$, namely for every definable subset
$X$ of $G$, $\mu(X)\geq 0$, $\mu(\emptyset) = 0$ and if $X, Y$ are
disjoint definable subsets of $G$ then $\mu(X\cup Y) = \mu(X) +
\mu(Y)$. (But note we do not require there be a finite bound on the
measures of definable sets).}
\end{Definition}

Note a difference with the usual situation: If
   $G$ is Ind-definable it may have NO type-definable subgroup of
bounded index (because $G$ itself is not type-definable). In any
case if $G$ has a smallest type-definable subgroup of bounded index
we will call it $G^{00}$ and say ``$G^{00}$ exists".

\vspace{5mm} \noindent As in Section 4, we will say that the
Ind-definable (over $\emptyset$) group $G =\bigcup_{n}G_{n}$ has
{\em finitely satisfiable generics} if there is a global complete
type $p(x)$ of $G$ (namely $p(x)\rightarrow ``x\in G_{n}"$ for some
$n$) such that every left translate of $p$ by an element of $G$ is
finitely satisfiable in some fixed small model $M_{0}$.

The material from Section 4 generalizes as follows:

\begin{Proposition}\label{6.3} Suppose the Ind-definable group $G$ 
has fsg. Then
\newline
(i) Any definable subset $X$ of $G$ is left generic iff right generic iff
every left (right) translate meets $M_{0}$.
\newline
(ii) there is a complete global generic type of $G$ (in fact $p$ as in
the definition of $fsg$ will be such, as well as any translate of $p$).
\newline
(iii) If $X$ is a definable subset of $G$ which is generic in $G$
and $X = X_{1} \cup X_{2}$ with $X_{i}$ definable, then $X_{1}$ or
$X_{2}$ is generic in $G$.
\newline
(iv) There is a smallest subgroup of $G$ which has bounded index and is
invariant over some small set of parameters.
\newline
(v) $G^{00}$ exists and
equals
$Stab(q)$ for each global generic type
$q$ of $G$. (Hence the cosets of $G^{00}$ in $G$ correspond to the
translates of $q$.)
\end{Proposition}
\pf Let $p$ be the type given by $fsg$. So for all sufficiently
large $n$, $``x\in G_{n}"\in p$. Likewise any definable subset $X$
of $p$ is in $G_{n}$ for sufficiently large $n$. So given a
definable left generic set $X$, there is an $n$ such that $``x\in
G_{n}"\in p$ and $X\subseteq G_{n}$. So (as finitely many left
translates of $X$ cover $G_{n}$) some left translate of $X$ is in
$p$ hence $X$ is in some left translate of $p$, so $X$ meets
$M_{0}$. Likewise every left translate of $X$ meets $M_{0}$. Now fix
$m$. Then for every $g\in G_{m}$, $gX$ meets $M_{0}$. By compactness
there are $g_{1},..,g_{k}\in G(M_{0})$ such that for every $g\in
G_{m}$, $gX$ contains one of the $g_{i}$. But then for every $g\in
G_{m}$, $g^{-1}\in Xg_{1}\cup..\cup Xg_{k}$. As $G_{m} =
G_{m}^{-1}$, we see that finitely many right translates of $X$ cover
$G_{m}$. As $m$ was arbitrary we conclude that $X$ is right generic.
The rest of (i) follows by the same argumentation (noting that every
right translate of $p^{-1}$ is finitely satisfiable in $M_{0}$).
\newline
(ii) follows from (i).
\newline
(iii) Again we may suppose that $X\subseteq G_{n}$ and $p(x)\models
``x\in G_{n}$. So some translate of $X$ is in $p$, so $X$ is in a
translate of $p$, so $X_{1}$ or $X_{2}$ is in the same translate of
$p$ so is generic.
\newline
At this point we see that the collection of nongeneric definable subsets
of $G$ is an ideal. Call this ideal ${\cal I}$.
\newline
(iv) Suppose $H$ to be a subgroup of $G$ of bounded index which is
$A$-invariant for some small set $A$. Then our given global generic type
$p$ determines a coset of $H$ in $G$ and every other coset of $H$ in $G$
corresponds to a translate of $p$. So the number of translates of $p$
bounds the index of $H$ in $G$. Hence there is smallest such $H$ (even
as $A$ varies).
\newline
(v) requires a little finesse. First let $X$ be any definable subset
of $G$. Let $Stab_{\cal I}(X)$ be as in Section 4, namely $\{g\in G:
gX\Delta X$ is nongeneric in $G$\}. So $Stab_{\cal I}(X)$ is a
subgroup of $G$, but on the face of it has no definability
properties. But we DO know that $Stab_{\cal I}(X)$ is invariant over
the parameters defining $X$, and also has bounded index in $G$ (as
generics meet $M_{0}$). Now fix a global generic type $q$. By what
we have just said, together with (iv), $\cap\{Stab_{\cal I}(X):X\in
q\} = H$ say is a subgroup of bounded index invariant over some
small set,  and $H$ is clearly contained in $Stab(q)$. But there is
an $n$ such that $q(x)\models ``x\in G_{n}"$, and therefore, as
above, $G_n$ and every $G_m$, $m>n$, are generic in $G$.

Clearly $Stab(q)$ is a subgroup of $G$ contained in $G_{n+1}$. Thus
$H\subseteq G_{n+1}$. Now for $X\in q$, let $Stab_{\cal I}^{n+1}(X)
= \{g\in G_{n+1}:gX\Delta X \mbox{is nongeneric}\}$. But this is
clearly type-definable (as we only have to say that finitely many
translates of $gX\Delta X$ do not cover $G_{n+2}$.) As $H =
\cap\{Stab_{\cal I}^{n+1}(X):X\in q\}$, it follows that $H$ is
type-definable. So we have constructed a type-definable subgroup of
$G$ of bounded index. By (iv) there is a smallest one, so $G^{00}$
exists. As in the earlier proof, $G^{00}$ must contain $Stab(q)$. So
$G^{00} = H = Stab(q)$.\qed

\vspace{5mm} \noindent

We can  easily generalize Proposition \ref{5.5} as well.
\begin{Proposition}\label{6.2} Assume that $T$ has the NIP, and $G = 
\bigcup_{n}G_{n}$
is an Ind-definable group (Ind-definable over $\emptyset$). Suppose
that $G$ HAS a type-definable subgroup of bounded index. Then it has
a smallest one, $G^{00}$ which is moreover normal and type-definable
over $\emptyset$.
\end{Proposition}
\pf Note by assumption (*) that if $H$ is a type-definable subgroup
of $G$, contained in $G_{n}$ say, then $G/H$ has bounded cardinality
iff $G_{n}/H$ does.

By our assumptions, without loss of generality there is a 
type-definable subgroup
$H$ of $G$ of bounded index, which is contained in $G_{0}$. The
proof of Proposition \ref{5.5} goes through word for word to give  a
type-definable subgroup $L_{0}$ of $G$ of bounded index which is
smallest among those contained in $G_{0}$. Likewise for each $n$
there is a type-definable subgroup $L_{n}$ of $G$ which is smallest
among those contained in $G_{n}$. But then $L_{n}\subseteq L_{0}$ so
$L_{n}$ is contained in $G_{0}$ so $L_{n} = L_{0}$. Thus $L_{0} =
G^{00}$. It is clearly normal and type-definable over
$\emptyset$.\qed

\begin{Lemma}\label{6.4} Assume that $G$ is an Ind-definable group as 
above and that $G^{00}$
is a minimal type-definable subgroup of bounded index. Let $\pi:G\to
G/G^{00}$ be the projection map and set $Y\subseteq G/G^{00}$ to be
closed if and only if for every $n$, $\pi^{-1}(Y)\cap G_n$ is
type-definable. Then these closed sets generate a locally compact
topology on $G$, making it into a topological group.

The compact sets in $G/G^{00}$ are those closed $Y$ such that
$\pi^{-1}(Y)$ is contained in $G_n$ for some $n$.
\end{Lemma}
\pf Left to the reader.\qed

\begin{Remark} In fact one can formulate the notion of a
``type-definable" equivalence relation $E$ on an Ind-definable set $Y$,
and assuming boundedly many classes one can define the ``logic topology"
on
$Y/E$ which will  be locally compact. As we will only require the 
group case as in 7.5, 
we leave details of the general case to the 
reader.
\end{Remark}

Finally we generalize Proposition \ref{5.6} to the Ind-definable
setting. Recall first that a left Haar measure on a locally compact
group $G$ is a left invariant Borel measure $\mu$ on $G$ such that
$\mu(X)$ is finite for $X$ compact and positive for $X$ open (so may
take value $\infty$ on some Borel sets). A left Haar measure exists
and is unique up to multiplication by a positive real.

\begin{Proposition}\label{6.6} Let $G$ be an Ind-definable group with finitely
satisfiable generics. Assume $T$ has NIP. Then there is a left invariant
Keisler measure on $G$ which is moreover finitely satisfiable in some
small model.
\end{Proposition}
\pf As in the proof of \ref{5.6}, we may assume $T$ to be countable.
Let $\textbf m$ be a right Haar measure on the locally compact group
$G/G^{00}$. Let $p(x)$ be a global generic type extending $G^{00}$.
Without loss of generality $G^{00}$ is contained in $G_{0}$.
\newline
We would like (as in the proof of Proposition \ref{5.6}) to define a
left invariant Keisler measure $\mu_{p}$ on $G$ by stipulating that
for any definable subset $X$ of $G$, $\mu_{p}(X) = {\textbf m}
(\{g/G^{00}: gX \in p\})$.

So fix a definable subset $X$ of $G$. Assume $X\subseteq G_{n}$. As
before, whether or not $gX$ is in $p$ depends only on $g/G^{00}$. So the
main point is to see that
$\{g/G^{00}:gX\in p\}$ is Borel and has finite
${\cal M}$-measure.

Note that if $gX\in p$ then $g\in G_{n+1}$ (as $g\in G_{0}\cdot 
G_{n}$) and so $gX\subseteq
G_{n+2}$. We copy the proof of Proposition \ref{5.6} but defining
now $U$ to be $\{Y\cap G_{n+2}(M_{0}):Y\in p\}$, and concluding that
$\{g/G^{00}:gX\in p\}$ is a Borel subset of the compact set
$G_{n+1}/G^{00}$ hence has finite ${\textbf m}$-measure. So we can
define $\mu_{p}$. Left invariance, finite additivity, and finite
satisfiability in $M_{0}$ are proved as before.\qed

\vspace{5mm}
\noindent

We conclude this interlude with a result that appears at first sight close
to the conjectures for compact groups, mentioned in the introduction.

\begin{Proposition} \label{7.8}  Let ${\bar M}$ be a
saturated $o$-minimal structure (expansion of a real closed field)
and $G$ a definably connected group definable in ${\bar M}$. Then:
\newline
There is a definably compact neighbourhood of the identity $G_{0} =
G_{0}^{-1}$ such that putting $G_{n} = G_{0}\cdot..\cdot G_{0}$, and
$G_{\infty} = \bigcup_{n}G_{n}$, then the Ind-definable group
$G_{\infty}$ has a unique smallest type-definable subgroup of
bounded index $G_{\infty}^{00}$ and the quotient $L= G_{\infty} 
/G_{\infty}^{00}$
  with the ``logic
topology" is a connected Lie group of the same dimension as the
$o$-minimal dimension of $G$.   \end{Proposition}

\proof We can identify some   neighborhood of $1$ in  $G$
with a neighborhood of 0 in $R^n$; write $*$ for
multiplication in $G$.  The only possible linear approximation to 
$x*y$ is $x+y$, by associativity
and the existence of differentiable inverse.  So letting $|x| = max 
|x_i|$, for any $C >0$,  for
all sufficiently small $e>0$, if $|x| \leq e$ and $|y|\leq e$ then
\begin{equation} \label{est}  |x * y - (x+y) | \leq C |(x,y)| \end{equation}
Take $C$ infinitesimal,  and then $e$ infinitesimal compared to it, 
and let $U  = \{x: |x| \leq e \}$,
$H = \{x: |x| \leq (1/n)e, n=1,2,...\}$.  Then by (\ref{est}) it is 
clear that $H$ is a type-definable
normal subgroup.   Let $G_0 = U \union U^{-1}$ in the sense of (*), 
so as to have it symmetric;
  Let $G_\infty$ be the Ind-definable group generated by $U$, or 
equivalently by $G_0$.
Modulo $H$, $*$ agrees with $+$ on $U$, indeed on $G_\infty$.   In 
particular $G_\infty/H \simeq \Rr^n$.

It remains only to show that $H$ equals $G^{00}$ precisely, i.e.  that
$G_{\infty} /G_{\infty}^{00}$
cannot have dimension {\em bigger} than  the
$o$-minimal dimension of $G$.  We postpone this to \S 10, see 
Corollary \ref{10.10}. \qed

\vspace{5mm}
\noindent
However,  note that the locally compact quotient we obtained is 
abelian; it is indeed a locally compact manifestation of the Lie 
algebra of $G$.  We feel that the canonical compact quotient of a 
definably compact group $K$ reflects better the structure of $K$; for 
instance $K/K^{00}$
is non-abelian if $K$ is non-abelian.  In the general case too, there 
should also be
a locally compact quotient whose structure is close to that of $G$.
We do not at the moment have a precise statement
of this, either in the compact or in the locally compact cases.

Note that the adjoint action $G \times L \to L$ is definable, in the 
sense of \S 2.
%

\section{Proof of the o-minimal conjectures}
We now use some of the preceding results to complete the proof of
the conjectures on definably compact definable groups in $o$-minimal
structures from \cite{Pillay}. In fact we will prove a bit more,
namely that such groups have $fsg$ and therefore, by \ref{5.6} are
definably amenable. Our main result (stated in the language of
   Definition \ref{2.1}) is:
\begin{Theorem} \label{7.1} Let ${\bar M}$ be a saturated $o$-minimal
expansion of a real closed field. Let $G$ be a definably
connected definably compact group definable in ${\bar M}$.
Then
\newline
(i) $G$ has $fsg$.
\newline
(ii) There is a definable surjective homomorphism
$\pi:G\rightarrow H$ from $G$ to a compact Lie group $H$
such that the Lie group dimension of $H$ equals the
$o$-minimal dimension of $G$, and moreover such that any
definable homomorphism from $G$ to a compact group factors
through $\pi$.
\end{Theorem}

Of course the $H$ in part (ii) of the theorem is precisely
$G/G^{00}$ equipped with the logic topology. We know from
\cite{BOPP} that $G^{00}$ exists and $G/G^{00}$ is, as a topological
group,  a compact connected Lie group. As discussed in
\cite{Peterzil-Pillay} we may assume that $G$ is a definable closed
subset of some ${\bar M}^{n}$ and that the group operation on $G$ is
continuous with respect to the induced topology on $G$.

We will prove Theorem \ref{7.1} by proving it  when $G$ is
commutative and when $G$ is ``semisimple", and then use Proposition
\ref{4.5} among other things to conclude the general case. For the
rest of this section $\bar M$ is a saturated $o$-minimal expansion
of a real closed field.

\begin{Lemma} \label{7.2} Theorem \ref{7.1} is true when $G$ is commutative.
\end{Lemma}
\pf We use additive notation for $G$. We first prove (ii). $T$ being
$o$-minimal has NIP. Also as $G$ is commutative it is amenable  so
in particular definably amenable. Also by \cite{Peterzil-Pillay} the
family of nongeneric definable subsets of $G$ forms an ideal ${\cal
I}$. We can apply Proposition \ref{5.7} to conclude that $Stab_{\cal
I}(X)$ is a type-definable subgroup of $G$ of bounded index for any
definable subset $X$ of $G$. It is explained in
\cite{Peterzil-Pillay} how this implies (ii), but
   we briefly recall the argument.
For each
$n$, we can find a definable subset $X_{n}$ of $G$ such that
the sets $X_{n}, X_{n}+c_{1},..,X_{n}+c_{r}$ form a partition
of
$G$, where $0,c_{1},..,c_{r}$ are the elements of order $n$
in
$G$. Then $Stab_{ng}(X_{n})$ contains no $n$-torsion (except
$0$). So if we know that each $Stab_{ng}(X)$ has bounded
index it will follow that $G^{00}$ is contained in every
$Stab_{ng}(X_{n})$, hence has no torsion. As $G^{00}$ is
divisible (see \cite{BOPP}), it follows that $G$ and
$G/G^{00}$ have isomorphic torsion.

  By a theorem of Edmundo and Otero (see \cite{Edmundo-Otero}),
the torsion of $G$ is isomorphic to the torsion of
$(S^{1})^{dim(G)}$. Hence the compact commutative Lie group
$G/G^{00}$ must also be $(S^{1})^{dim(G)}$. So (ii) is proved.

\vspace{2mm} \noindent Now for (i). Let $\sim_{\cal I}$ be the
equivalence relation: ``$X\Delta Y$ is nongeneric" on definable
subsets of $G$. By Proposition \ref{5.7}(i) there are only boundedly
many definable subsets of $G$ up to $\sim_{I}$. (Note this already
proves that $G$ has a bounded number of generic types.) Thus there
is a small model $M_{0}$ such that  $G$ is defined
over $M_{0}$
and
for every generic definable
subset $Y$ of $G$ there is an $M_{0}$-definable subset $X$ of $G$
such that $Y\sim_{I}X$.  To prove that $G$ has $fsg$ it is clearly enough
(given the existence of generic types) to prove that every generic
definable subset $Y$ of $G$ meets $G(M_{0})$. So let $Y\subseteq G$
be definable and generic. \vspace{2mm}

   \noindent{\em Claim 1.}
There exists a definable subset $Y'\subseteq Y$ which is closed
(in $G$, so in ${\bar M}^{n}$) and still generic.

\pf First, we may replace $Y$ by its interior. Now, for every
$\epsilon
>0$ we consider the set $Y_{\epsilon}$ of all $y\in Y$ whose distance from the
frontier of $Y$ is greater than $\epsilon$ (in the sense of $\bar
M^n$). Because the frontier of $Y$ is not generic, there is some
$\epsilon>0$ for which $Y_{\epsilon}$ is generic, and we take it to
be $Y'$.

\vspace{2mm}
\noindent
  So we may assume $Y$ to be closed. Let
$X$ be an $M_{0}$-definable subset of $G$ such that $Y\sim_{ng}X$.
We may clearly assume $X$ to be closed (as $cl(X)\setminus X$ is
nongeneric). Hence $X\cap Y$ is closed. Let $Z = X\setminus Y$.
\vspace{2mm}

  \noindent{\em Claim 2.} The set of $M_{0}$-conjugates
of $X\cap Y$ is finitely consistent.
\newline
{\em Proof.} Otherwise (as $X$ is $M_{0}$-definable) finitely
many
$M_{0}$-conjugates of
$Z$ cover $X$. But $Z$ is nongeneric in $G$ as is any
$M_{0}$-conjugate of $Z$. So $X$ is the union of finitely
many nongenerics, while itself being generic. This is a
contradiction.

\vspace{2mm} \noindent By Claim 2 and Theorem 2.1 of
\cite{Peterzil-Pillay} (which comes out of Dolich's work
\cite{Dolich}), $X\cap Y$ meets $M_{0}$, as does $Y$. This completes
the proof of (i) and of Lemma \ref{7.2}.\qed

\vspace{5mm} \noindent Let $G$ be definable in ${\bar M}$. We will
say that $G$ has {\em very good reduction} if it is definably
isomorphic, in $\bar M$, to a group $G_1$ with the following
property: There is a sublanguage $L_{0}$ of the language $L$ of
${\bar M}$ which contains $+,\cdot$ and there is an elementary
substructure $M_{0}$ of ${\bar M}|L_{0}$ whose underlying set is
$\R$, and such that $G_1$ is definable by an $L_{0}$-formula with
parameters from $M_{0}$, i.e. from $\R$.  (But note that $M_{0}$
need not be expandable to an elementary substructure of ${\bar M}$.)

\noindent{\bf Remark}  This notion of very good reduction is related
to, but not identical with the algebraic-geometric notion in the
case of saturated real closed fields and the natural valuation. In
any case it is important to note that even if $R$ is a saturated
real closed field, there will be definable groups, even real
algebraic ones which do not have very good reduction in the model
theoretic sense above. Indeed, as was shown in
\cite{Peterzil-Starchenko2}, if $R$ is a sufficiently saturated real
closed field then not all elliptic curves over $acl(R)$ are
definably isomorphic to each other (as groups). In fact this remains
true even in an  expansion of $R$ to a structure $R_{an}$
elementarily equivalent to ${\mathbb R}_{an}$. Now, in ${\mathbb
R}_{an}$ all definable compact abelian groups of fixed dimension
(defined over $\mathbb R$) are definably isomorphic to each other,
therefore, even in $R_{an}$ not all elliptic curves over $acl(R)$
have very good reduction.

\begin{Lemma} \label{7.3} Theorem \ref{7.1} holds when $G$ has very 
good reduction.
\end{Lemma}
\pf Part (ii) of the theorem is precisely Fact 4.1 of
\cite{Peterzil-Pillay}. The fact that the nongeneric sets form an
ideal was proved in \cite{Peterzil-Pillay}, but this as well as the
rest of (i) follows directly from Proposition 4.6 in the same paper,
(which itself depends on results of Berarducci and Otero \cite{BO}).
More precisely (with above notation) 4.6 of \cite{Peterzil-Pillay}
states among other things that if $X\subset G$ is definable (in
${\bar M}$) then $X$ is left generic iff right generic iff $X$
contains an open set which is $L_{0}$-definable over $M_{0}$. This
on the one hand implies that there exists a complete generic type,
and on the other hand that if we pick $M_{1}$ to be any elementary
substructure of $\bar M$ which contains $M_{0}$ then any generic
definable subset of $X$ meets $M_{1}$. Thus $G$ has $fsg$.\qed

\vspace{5mm} \noindent {\em PROOF OF THEOREM \ref{7.1}.}
\newline
Let $G$ be an arbitrary definable, definably connected, definably
compact group in ${\bar M}$. We prove the theorem by induction on
$dim(G)$. If $G$ is ``semisimple", namely has no proper connected
infinite definable normal commutative subgroup, then by \cite{PPS1},
$G$ is an almost direct product of finitely many  almost definably
simple groups $G_{1},..,G_{k}$. (``Almost definably simple'' means
that the group is noncommutative and the quotient by some finite
normal subgroup is definably simple.) Now by \cite{PPS2} (see the
proof of $(2) \Rightarrow (3)$ in Theorem 5.1 there),  any definably
simple group is definably isomorphic to some semialgebraic group
defined over $\R$. In particular, a definably simple group has very
good reduction. It easily follows from Lemma \ref{7.3} that Theorem
\ref{7.1} holds for a semisimple $G$.

Thus we may assume that $G$ has an infinite, definably  connected
normal commutative subgroup $N$. By \ref{7.2}, the theorem is true
of $N$, so we may assume $N\neq G$. By induction,  the theorem is
true for $G/N$, so by Proposition \ref{4.5}, $G$ has $fsg$.

All that is left to do is to prove that the dimension of the compact
Lie group $G/G^{00}$ equals the $o$-minimal dimension of $G$. Notice
first that the image of $G^{00}$ under the projection onto $G/N$ is
necessarily $(G/N)^{00}$ (on one hand this image contains
$(G/N)^{00}$; on the other hand the pre-image of $(G/N)^{00}$ is of
bounded index and therefore contains $G^{00}$). Thus, it suffices to
show that $G^{00}\cap N = N^{00}$. By \cite{BOPP} it is enough to
prove: \vspace{2mm}

\noindent{\em Claim.} $G^{00}\cap N$ is torsion-free.

\pf Fix $n$. Let us first choose a definable subset $X$ of $N$ such
that $N$ is the disjoint union of the translates of $X$ by the
distinct $n$-torsion points $1,g_{1},..,g_{r}$ say of $N$. (As usual
$X$ is obtained by considering the surjective endomorphism $\pi:
x\to nx$ of $N$ with itself, which has finite kernel, and use the
existence of definable Skolem functions.) Likewise, using definable
Skolem functions, we can find a definable subset $D$ of $G$ which
meets every coset of $N$ in $G$ in a unique point. It follows that
the definable sets $XD$, $g_{1}XD$,..,$g_{r}XD$ are disjoint and
cover $G$. By Corollary \ref{4.3}, $G^{00}$ is contained in
$Stab_{\cal I}(X)$, and clearly the latter does not contain any of
$g_{1},..,g_{r}$. As $n$ was arbitrary, it follows that $G^{00}\cap
N$ is torsion-free. This completes the proof of Theorem
\ref{7.1}.\qed

\vspace{5mm} \noindent By Proposition \ref{5.6} we conclude also:
\begin{Corollary} Let $G$ be a definably compact group
definable in $\bar M$. Then $G$ is definably amenable. In
fact there is a left invariant Keisler measure on $G$
which is finitely satisfiable in some small model.
\end{Corollary}

\noindent{\bf Remarks}

\noindent 1. In the very last step of the above  proof we showed
that, under the given assumptions, $G^{00}\cap N=N^{00}$. This is
not true in general, even if we assume that $G$, $G/N$ and $N$ all
have  NIP and fsg. Indeed, consider the group $G=\la {\mathbb
C},+\ra \oplus S^1$ ($S^1$ the circle group), with predicates for
$S^1$ and all its semialgebraic subsets (but not for $\mathbb C$!).
We have $G^{00}=G$, but $(S^1)^{00}$ is nontrivial.

\noindent 2. Our proof of Theorem \ref{7.1} depends in a crucial
manner on the result \cite{Edmundo-Otero} describing the torsion in
definably compact commutative groups, which itself relies on quite
intricate tools from algebraic topology. It would be desirable to
have a ``direct" proof of the latter in the spirit of the current
paper. In fact we do have a reasonably elementary proof of the {\em
existence} of torsion points (in commutative definably compact
definably connected groups), which we sketch here:

\noindent (i) Using definable compactness, find a definable
$X\subset G$ such that both $X$ and its complement $X^{c}$ are
generic (this can be done similarly to the proof of Claim 1 above),
\newline
(ii) It follows that $Stab_{ng}(X) \neq G$, and thus (as we saw that
$Stab_{ng}(X)$ has bounded index), $G^{00} \neq G$.
\newline
(iii) Since $G/G^{00}$ is a compact connected commutative nontrivial
Lie group (\cite{BOPP}) it has torsion, and since $G^{00}$ is
divisible (\cite{BOPP}), $G$ itself has torsion. \vspace{2mm}

\noindent 3. Notice that if a definable $G$ in an o-minimal
structure has fsg then it necessarily implies that $G$ is definably
compact. Indeed, if $G$ were not definably compact then by
\cite{Peterzil-Steinhorn}, $G$ has a definable one dimensional,
ordered subgroup $H$. Let $D\subseteq G$ be a definable set
containing one representative for each coset of $H$, and let
$I=(0,\infty)\subseteq H$. Then $D\cdot I$ is nongeneric in $G$ and
so is its complement, contradicting \ref{4.2}.

\noindent 4. The proof of the o-minimal group conjecture that we
give here depends in the ambient real closed field in two different
ways. Firstly, in order to ensure that our group can be embedded as
a topological group into some $R^n$ (see a discussion in
\cite{Peterzil-Pillay}. Secondly (and more substentially) the above
count of torsion points, By Edmundo and Otero  was only carried out
for expansions of real closed fields. The conjecture was proved
separately for groups definable in ordered vector spaces over
division rings (see \cite{Onshuus}. \cite{Elef}).

\section{Compact domination}
The third author has mentioned in previous papers that the
$o$-minimal conjectures (solved in the last section) have the heuristic
content that the map $G\to G/G^{00}$ should be a kind of intrinsic
``standard part map". It is reasonable to attempt to give some concrete
mathematical meaning to this, namely to come up with a model theory of
``standard-part-like" maps (in a tame context). So we introduce the
notion ``compact domination".   It is analogous to ``stable domination"
from \cite{Hrushovski} which was introduced
with algebraically closed valued fields as a central example.
We relate compact domination to the
existence and uniqueness (and smoothness) of suitable Keisler measures,
and prove that in the cases we understand well (very good reduction and
dimension $1$) definably compact groups in $o$-minimal structures {\em
are} compactly dominated (by $G/G^{00}$).

We begin by working in a saturated model ${\bar M}$ of an arbitrary
theory. When we say compact we mean compact Hausdorff. $G$ denotes a
definable (or even type-definable) group. We use freely the notion
from Section 2 of a {\em definable map} from $X$ to a compact space.


\begin{Definition}\label{8.1}{\em (i) Suppose $X$ is
  type-definable, $\pi:X\to C$ is a definable surjective map from $X$ 
to a compact space $C$,
and $\mu$ is a probability measure on $C$. We
say that {\em $X$ is compactly dominated by $(C,\mu,\pi)$} if for any 
definable 
(that is relatively definable with
parameters) subset $Y$ of $X$, and for every $c\in C$ outside a set
of $\mu$ measure zero, either $\pi^{-1}(c)\subseteq Y$ or
$\pi^{-1}(c)\subseteq X\setminus Y$. Namely,
  $$\mu(\{c\in C:\pi^{-1}(c)\cap Y \neq\emptyset \mbox{ and }
\pi^{-1}(c)\cap (X\setminus Y) \neq \emptyset\}) = 0.$$ 
\newline
(ii) Let $G$ be a type-definable  group. We say
that {\em $G$ is compactly dominated as a group},
if $G$ is compactly dominated over  by $(H,{\textbf
m},\pi)$ where $H$ is a compact group, $\textbf m$ is the unique
normalized Haar measure on $H$ and $\pi$ is a group homomorphism.}
\end{Definition}

Note that in (i) above the set $\{c\in C:\pi^{-1}(c)\cap Y \neq
\emptyset$ and $\pi^{-1}(c)\cap (X\setminus Y)\neq\emptyset\}$ is a
closed subset of $C$, hence measurable.

When we work with a definable group $G$, we always refer to compact
domination in the group sense.

\begin{Question} To what extent does the
definition of compact domination depend on the choice of ``measure zero'' as
the  notion of ``smallness" in $C$? \end{Question}

  It would be interesting to investigate
   other possibilities.   Smallness notions based on Baire category or 
dimension are more natural
  since they depend only on the topology; but in the context of groups 
the Haar measure
  also depends only on the topology and group structure, and connects naturally
  to the topics discussed in this paper.
   It would be nice if for groups these notions turned out to be equivalent.

Let  $P$ be compactly dominated via  $\pi:P\rightarrow
C$, where $P$ and $\pi$ are (type-) defined over $\emptyset$.   We 
will say ``$\theta(x,b)$ holds for almost all $x \in P$''
if $\mu(\pi( \{x: \neg \theta(x,b) \})=0$.
We can write:
$(d_Px)\theta(x,b)$ for this.  Note that this gives an partial type:
$\{b:  (d_Px)\theta(x,b) \}$ is type-definable over $\emptyset$. 
Indeed let $\{W_i\}_{i \in I}$ be
the set of all closed subsets of $C$ of positive measure; then $\pi 
^{-1} (W_i) = \cap_j W_{ij}$
for some definable sets $W_{ij}$.  Now
$\neg (d_Px)\theta(x,b)$  iff  $\mu(\pi(\theta(x,b)))>0$  iff
$\pi(\theta(x,b))$ contains a closed set $W_i$ of measure $>0$, iff 
for some $i,j$ $\pi(\neg \theta(x,b))$ contains $W_{ij}$.  The case 
of Baire category is similar.

This is again in analogy with the stably dominated case, where one 
obtains definable types.

One could ask to what extent $C$ is determined by $P$?  If $P$ is 
compactly dominated
via $\pi_i: P \to C_i$, there exist continous maps $f_1: C_1' \to C_2$
and $f_2: C_2' \to C_1$, where $C_i'$ is a large subset of $C_i$, such that
$f_1 \pi_1 = f_2$ for all $x \in \pi_1 ^{-1} (C_1')$, and dually. 
However in general
$f_1,f_2$ are not inverses of each other.

\begin{Proposition}\label{8.3} Suppose $G$ is compactly dominated by 
$(H,\pi)$. Then
\newline
(i) $G$ has finitely satisfiable generics, and
\newline
(ii) $G^{00}$ exists and equals $Ker(\pi)$.
\end{Proposition}
\pf Let us assume that $G$ is compactly dominated (over $\emptyset$
say) by the data. We go through various claims which eventually
yield (i) and (ii). $Y$ will denote a definable subset of $G$ and
$\pi'(Y) = \{h\in H:\pi^{-1}(h)\subseteq Y\}$.
\newline
{\em Claim 1.} $\pi'(Y)\subseteq
\pi(Y)$, $\pi(Y)$ is closed and $\pi'(Y)$ is open.
\newline
\pf  Clear.

\vspace{2mm} \noindent {\em Claim 2.} ${\textbf m}(\pi(Y)) =
{\textbf m}(\pi'(Y))$.
\newline
\pf Because by the definition of compact domination ${\textbf
m}(\pi(Y)\setminus\pi'(Y)) = 0$.

\vspace{2mm}
\noindent
{\em Claim 3.} The following are equivalent:
\newline
(a) $Y$ is left (right) generic,
\newline
(b) $\pi(Y)$ is left(right) generic,
\newline
(c) ${\textbf m}(\pi(Y)) > 0$,
\newline
(d) $\pi'(Y)$ is nonempty.

\vspace{2mm} \noindent {\em Proof.} (a) implies (b) implies (c) are
clear. Suppose now that (c) holds. Then by Claim 2, ${\textbf
m}(\pi'(Y))
>  0$, so in particular (d) holds. Now assume (d). Then $Y$ contains
a coset of $Ker \pi$, which is type-definable of bounded index, and
hence $Y$ is left and right generic.

\vspace{2mm} \noindent {\em Claim 4.} If $Y = Y_{1}\cup Y_{2}$
(where the $Y_{i}$ are definable) and $Y$ is generic, then $Y_{1}$
or $Y_{2}$ is generic.

\vspace{2mm}\pf  By Claim 3, ${\textbf m}(\pi'(Y))>0$, but the
compact domination assumption implies that ${\textbf
m}(\pi'(Y))={\textbf m}(\pi'(Y_1))+{\textbf m}(\pi'(Y_2))$, so again
by Claim 3 we are done.

\vspace{2mm} \noindent Let $M_{0}$ be an elementary substructure of
${\bar M}$ containing representatives of each coset of $G$ modulo
$Ker(\pi)$. \vspace{2mm}

\noindent {\em Claim 5.} If $Y$ is generic in $G$ then $Y$ meets
$M_{0}$.

\noindent \pf By Claim 3, $Y$ contains a whole coset of $Ker(\pi)$.

\vspace{2mm}
\noindent
By Claim 4  there is a global generic type $p$ of $G$. Every
translate of $p$ is also generic so by Claim 5 is finitely satisfiable in
$M_{0}$. Thus $G$ has $fsg$ giving part (i).

\vspace{2mm} \noindent Let ${\cal I}$ be the ideal of nongeneric
definable sets ( which exists by Claim 5.)
\newline
{\em Claim 6.} $Ker(\pi) \subseteq Stab_{\cal I}(Y)$.
\newline
{\em Proof.}  Let $g\in Ker(\pi)$. Then $\pi(Y\Delta gY)\subseteq
(\pi(Y)\setminus \pi'(Y))\cup (\pi(gY)\setminus \pi'(gY))$. By Claim 2
the latter has Haar measure $0$, hence by Claim 3, $Y\Delta gY$ is
nongeneric.

\vspace{2mm} \noindent By Corollary \ref{4.3}, $G^{00}$ exists and
equals the intersection of all $Stab_{\cal I}(Y)$. Since $ker \pi$
has bounded index in $G$, by Claim 6, $G^{00}$ equals
$Ker(\pi)$.\qed

\vspace{5mm}
\noindent
Note that it follows that if $G$ is $\emptyset$-definable and compactly
dominated over some parameters then it is compactly dominated over any
model (as $G^{00}$ is type-definable over $\emptyset$). We now aim towards
the appropriate analogue of ``existence and uniqueness of Haar measure"
for compactly dominated groups. We begin with a group-free version:

\begin{Proposition}\label{8.4} Let $X$ be type-definable over $\emptyset$, and
compactly dominated over $\emptyset$ by $(C,\mu,\pi)$. Then:
\newline
(i) There is a unique Keisler measure $\mu'$ on $X$ with the
property that $\mu(D) = \mu'(\pi^{-1}(D))$ for any closed
$D\subseteq C$.
\newline
(ii) The Keisler measure $\mu'$ from (i) is smooth (over $\emptyset$).
\end{Proposition}
\pf We first start with an explanation. Given a Keisler measure
$\nu$ on a definable or type-definable set, we can uniquely extend
$\nu$ to a countably additive measure on the $\sigma$-algebra whose
underlying ``closed" sets are the type-definable subsets of $X$.
(This was discussed and referenced in section 2.) So, as $\pi^{-1}(C)$ is
type-definable (over $\emptyset$), then in (i) $\mu'(\pi^{-1}(C))$
makes sense, for a Keisler measure $\mu'$. In fact it is precisely
the infimum of the $\mu'(Y)$ for $\emptyset$-definable $Y$
containing $\pi^{-1}(C)$.

\vspace{2mm}
In any case, let us first show the existence of $\mu'$: For $Y$ a
(relatively) definable subset of $X$, put $\mu'(Y) = \mu(\pi(Y))$. Note
that $\mu'$ DOES satisfy the condition in (i): for if $D\subseteq C$ is
closed, and $Y= \pi^{-1}(D)$, then $Y$ is type-definable so equals
$\cap_{i}Y_{i}$ where $Y_{i}$ are (relatively) definable subsets of $X$.
Let $D_{i} = \pi(Y_{i})$. Then $D_{i}$ is closed in $C$ and
$\cap_{i}D_{i} = D$. We may assume that the family $(Y_{i})_{i}$ is
closed under finite intersections. It follows that $\mu(D) =
inf\{\mu(D_{i}):i\in I\} = inf\{\mu'(Y_{i}):i\in I\} = \mu'(Y)$.

We
must check  finite additivity of $\mu'$. But if $Y_{1}, Y_{2}$ are
disjoint definable subsets of $X$, then (by compact domination)
$\mu(\pi(Y_{1})\cap \pi(Y_{2})) = 0$, hence $\mu'(Y_{1}\cup Y_{2}) =
\mu'(Y_{1}) + \mu'(Y_{2})$.

\vspace{2mm} Now for uniqueness: Suppose $\mu''$ is another Keisler
measure on $X$ such that $\mu(D) = \mu''(\pi^{-1}(D))$ for any
closed $D\subseteq C$. Let $Y$ be an arbitrary definable subset of
$X$. Then, since $\pi^{-1}\pi'(Y)\subseteq Y\subseteq
\pi^{-1}\pi(Y)$, we have
$$\mu(\pi'(Y))=\mu''(\pi^{-1}\pi'(Y))\leq \mu''(Y)\leq
\mu''(\pi^{-1}\pi(Y))=\mu(\pi(Y)).$$ But $\mu(\pi'(Y))=\mu(\pi(Y))$,
hence $\mu'(Y)=\mu''(Y)$. So we have proved (i).

Recall that the smoothness of $\mu'$ over $\emptyset$ means by
definition that $\mu'|\emptyset$ has precisely one extension to a
Keisler measure on (all definable subsets of) $X$. However, since
$\mu'|\emptyset$ satisfies the assumptions of (i) it follows that it
has a unique extension.\qed

\begin{Theorem}\label{8.5} Suppose $G$ is compactly dominated. Then $G$ has a
unique left invariant Keisler measure, which is moreover right
invariant and smooth.
\end{Theorem}
\pf Let $\pi:G\rightarrow H = G/G^{00}$. As before ${\textbf m}$
denotes the Haar measure on $H$.

Let $\mu'$ be as in Proposition \ref{8.4} and its proof, namely for
definable $X\subseteq G$, $\mu'(X)$ is by definition ${\textbf
m}(\pi(X))$. Note that $\mu'$ will be both left and right invariant,
as ${\textbf m}$ is. By Proposition \ref{8.4} $\mu'$ is also smooth.

Now suppose $\mu''$ is another left invariant Keisler measure on
$G$. Let $M_{0}$ be a model over which $\pi$ is definable. By
\cite{Keisler}, $\mu''|M_{0}$ extends uniquely to a countably
additive measure on the $\sigma$-algebra of subsets of $G$ generated
by the $M_{0}$-type-definable sets. We still call this $\mu''|M_{0}$
and note it is left invariant. But then $\mu''|M_{0}$ induces a left
invariant countably additive measure on $H$: namely for $B$ a Borel
subset of $H$, define its measure to be $\mu''(\pi^{-1}(B))$. By
uniqueness of Haar measure, this latter measure has to agree with
${\textbf m}$. Hence we have shown that ${\textbf m}(C) =
\mu''(\pi^{-1}(C))$. By Proposition \ref{8.4} (i), $\mu'' = \mu'$.
This completes the proof.\qed

\section{$o$-minimality and compact domination}

Let ${\bar M}$ denote now a saturated $o$-minimal expansion of an
ordered divisible group $R$.

Beraducci and Otero, in their paper \cite{BO}, prove in effect, (for
$o$-minimal expansions of real closed fields) that the unit $n$-cube
$I^{n}$ in ${\bar M}$ is compactly dominated, with respect to the
standard part map to $I^{n}(\R)$ equipped with Lebesgue measure.
This is not stated explicitly in their paper, but follows from it.
In any case we give below another proof of this fact (omitting the
real closed field assumption), using the following  beautiful
theorem of Baisalov and Poizat (Recall that {\em a weakly
$o$-minimal structure} is an ordered structure in which every
definable subset of the linear ordering is a finite union of convex
sets): \vspace{2mm}

  \noindent{\bf
Theorem}(\cite{BP}) {\em If  the saturated $o$-minimal structure 

$\bar M$ is expanded by any number of
convex subsets of $\bar M$ then the resulting structure is weakly
$o$-minimal}. \vspace{1mm}

\noindent  Some notation: We let $\R$ denote a fixed copy of
  the reals, which we may assume
  is a subgroup of $R$ (in particular, we have a copy of $\mathbb Q$
  in $R$). Let $Fin$ denote the set of finite elements of $R$ (i.e.
absolute value less than $n$ for some $n\in {\N}$) and $Inf$ the set
of infinitesimals of $R$ (absolute value $< 1/n$ for all $n\in
{\N}$). Let $\pi$ denote the ``standard part map'' from $Fin$ onto
$Fin/Inf$. Since $Fin/Inf$ is archimedean (and $\bar M$ saturated)
we can identify $Fin/Inf$ with $\R$.

Let $\la {\bar M},Fin, Inf\ra$ be the structure ${\bar M}$ equipped
with unary predicates for $Fin$ and $Inf$. Then the quotient group
$Fin/Inf$ is interpretable in it, and  $\pi$ induces a canonical
bijection $i: Fin/Inf \rightarrow \R$.
\begin{Definition} By ${\R}_{ind}$ (standing for ``$\R$ with the induced
structure") we mean the structure whose universe is $\R$ and whose
relations are precisely the images under $i$ of subsets of
$(Fin/Inf)^{n}$ which are definable (with parameters) in $({\bar
M},Fin, Inf)$.
\end{Definition}

\begin{Lemma}\label{9.2} ${\R}_{ind}$ is $o$-minimal (in fact is an $o$-minimal
expansion of the ordered group of $\R$).
\end{Lemma}
\pf It is clear that $<$ and the graphs of $+$ and $\cdot$ are among
the basic relations on ${\R}_{ind}$.

By \cite{BP} the structure $\la {\bar M}, Fin, Inf\ra$ is weakly
$o$-minimal. Let $X\subseteq {\R}$ be definable in ${\R}_{ind}$.
Then clearly $\pi^{-1}(X)$ is definable in $\la {\bar M},Fin,
Inf\ra$, so is a finite union of convex sets. So $X$ has finitely
many connected components. Thus ${\R}_{ind}$ is $o$-minimal.\qed

\begin{Lemma}\label{9.3} Let $X \subset Fin^{n}$ be definable in 
${\bar M}$ with
$dim(X) < n$. Then $dim(\pi(X))<n$ (in the o-minimal structure
${\R}_{ind}$).
\end{Lemma}
\pf The proof is by induction on $n$, and is immediate for $n=1$.
For an arbitrary $n$, we may assume by cell decomposition that $X$
is the graph of a continuous definable function $f:C\to R$, where
$C$ is a definable open set in $R^{n-1}$. By o-minimality of
$\R_{ind}$, if $dim(\pi(X))=n$  then it must contain the closure of
a subset $U\times (q_1,q_2)$, for $U$ an open rectangular box of
rational coordinates (which we may assume is contained in $C$) and
$q_1,q_2\in \mathbb Q$.

Consider an arbitrary $x\in U(R)$ and $r$ a rational number in
$(q_1,q_2)$. By assumptions, there exist $x_1, x_2$ infinitesimally
close to $x$ such that $f(x_1), f(x_2)$ are infinitesimally close to
$q_1,q_2$, respectively. But then, by continuity, there exists an
$x'$ infinitesimally close to $x$ such that $f(x')=r$. It follows
that $\pi(\{x\in U(R):f(x)=r\})=U$, which by induction implies that
the set  $\{x\in U(R):f(x)=r\}$ has an interior in $R^{n-1}$. This
can be done for any rational $r\in(q_1,q_2)$, contradiction.\qed

\begin{Theorem}\label{9.4} Let $I^{n}$ be the unit $n$-cube in $R^{n}$, $\pi$
the standard part map from $I^{n}$ to $I^{n}(\R)$, and $\mu$ the Lebesgue
measure on $I^{n}(\R)$. Then
$I^{n}$ is compactly dominated (in ${\bar M}$) by $(I^{n}(\R),
\mu, \pi)$.
\end{Theorem}
\pf Let $X\subseteq I^{n}$ be definable in ${\bar M}$. Let $Y$ be
the frontier of $X$ (the set of $x$ such that every neighbourhood of
$x$ contains points both in $X$ and not in $X$). Then $dim(Y) < n$.
So $dim(\pi(Y)) < n$ by Lemma 10.3. As $\pi(Y)$ is definable in the
$o$-minimal structure ${\R}_{ind}$, it follows that the Lebesgue
measure of $\pi(Y)$ is $0$. Note also that $\pi(Y)$ is closed. For
$c\in I^{n}(\R)$, the type-definable set $\pi^{-1}(c)$ is {\em
definably connected} (cannot be written as the union of two
relatively open relatively definable subsets). So for $c\in
I^{n}(\R) \setminus \pi(Y)$, either $\pi^{-1}(c)$ is contained in
$X$ or contained in the complement of $X$. This proves compact
domination.\qed

${\bar M}$.

We are now in a position to state a rather finer version of the
conjectures from \cite{Pillay}.
  As before $\pi$ denotes the homomorphism
from $G$ onto $G/G^{00}$ and ${\textbf m}$ denotes Haar measure on
$G/G^{00}$.
\newline

\noindent{\em Compact Domination Conjecture.} Any definably compact
group $G$ (definable in a saturated $o$-minimal expansion of a real
closed field) is compactly dominated (by the compact Lie group
$G/G^{00}$, with its Haar measure $\textbf m$).

\vspace{5mm} \noindent
Note that, by 9.3, if $G$ (definably compact 
in saturated $o$-minimal expansion of a real closed field) is 
compactly dominated by $H$, then $H$ has to coincide with the compact 
Lie group $G/G^{00}$.

\vspace{2mm}
\noindent
  The following lemma allows us to reduce the Compact Domination Conjecture
to a simpler statement.

\begin{Lemma} \label{9.6} Suppose $G$ is definably compact with 
$dim(G) = n$, and
suppose that whenever $Y\subseteq G$ is definable
   and $dim(Y) < n$, then ${\textbf m}(\pi(Y)) = 0$. Then $G$ is compactly
dominated by $G/G^{00}$.
\end{Lemma}
\pf Note that $G$ here is equipped with its ``definable topology".
We make use of a key result from \cite{BOPP} which says that
$G^{00}$, and each translate of it, are definably connected. It
follows that if $X\subseteq G$ is definable, and $Y$ is the frontier
of $X$ in $G$ (which has dimension $< n$) then for all $c\notin
\pi(Y)$, $\pi^{-1}(c)$ is either contained in $X$ or disjoint from
$X$. Now, just like in the proof of \ref{9.4}, we obtain compact
domination.\qed

  The above conjecture, if proven true, will resolve an
intriguing open problem regarding the connection between generic
sets and torsion points.
\begin{Proposition} Assume that $G$ a definable abelian group
in $\bar M$ and that $G$ is compactly dominated by $G/G^{00}$ (with
its Haar measure). Then every definable generic subset of $G$
contains a torsion point. In particular, if $X\subseteq G$ is
generic then there are finitely many torsion points $g_1,\ldots,
g_k$ such that $G=\bigcup_i g_iX$.
\end{Proposition}
\pf If $X\subseteq G$ is generic then, by Claims 1 and  3 in the
proof of \ref{8.3}, $\pi'(X)=\{g/G^{00}: gG^{00}\subseteq X\}$ is
open in $G/G^{00}$ and therefore contains a torsion point. Since
$G^{00}$ is divisible and torsion-free, the coset $gG^{00}$, and
therefore $X$, contain a torsion point. The rest easily follows.\qed
\vspace{2mm}

  There is very little we currently know about the
consequences of the above proposition. Indeed, we don't even know
that every large set (namely, the complement of a definable subset
of $G$ of small dimension) contains a torsion point.

\begin{Theorem}\label{9.7} Let $G$ be a definably compact group definable in
an $o$-minimal  ${\bar M}$. Then $G$ is compactly dominated in
either of the cases
\newline
(i) $\bar M$ expands a real closed field and $G$ has very good
reduction.
\newline
(ii) $dim(G) = 1$.
\end{Theorem}
\pf Case (i): We assume that there is a sublanguage $L_{0}$ of $L$
such that $G$ is defined in $L_{0}$ over the elementary substructure
$M_{0} = \la \R,+,<,..\ra$ of ${\bar M}|L_{0}$. Assume $dim(G) = n$.
Then $G$ has a covering by finitely many charts $U_{1},..,U_{r}$,
each of which is definably homeomorphic via some $f_{i}$ to an open
definable subset $V_{i}$ of $I^{n}$ (all definable in $L_{0}$ over
$M_{0}$). Let ${\R}_{ind}$ be as above. As was pointed out earlier,
$G^{00}$ is exactly the collection of all elements in $G$ that are
infinitesimally close to $e$. Thus we identify $G/G^{00}$ with
$G({\R}_{ind})$. Suppose $Y\subseteq G$ is definable with $dim(Y) <
n$. Then working in the charts and using \ref{9.3} we see that
$dim(\pi(Y)) < n$ in the $o$-minimal structure ${\R}_{ind}$. Then
clearly ${\textbf m}(\pi(Y)) = 0$. (For example, working in the
charts the Lebesgue measure of $\pi(Y)) = 0$, so the Haar measure
must be $0$ too.) Now apply \ref{9.6}.

Case (ii). If $dim(G) = 1$ then any definable subset $Y$ of $G$ of
dimension $< 1$ is finite, so $\pi(Y)$ is finite too hence has Haar
measure $0$. Again apply \ref{9.6}.\qed

\begin{Corollary} Suppose $G$ is as in Theorem \ref{9.7} Then there is a
unique invariant Keisler measure on $G$, which is moreover smooth.
\end{Corollary}
\pf By \ref{8.5} and \ref{9.7}.\qed

\vspace{5mm}
\noindent
Finally we return to the promised completion of the proof  of 
Proposition \ref{7.8}, this time as an
illustration of the compact domination conjecture.
Actually the dominating group is {\em locally compact}
in this case; the modification of the definition is evident.  We show initially
that $G_\infty$ is  (locally) compactly dominated via $G_\infty \to 
G_\infty/H$; as a
bi-product, this gives $H=G_\infty^{00}$.

\begin{Proposition} \label{10.9}  Let $G,H$ be as in Proposition 
\ref{7.8}.  Then $G_\infty$
is (locally) compactly dominated via $G_\infty \to G_\infty/H$. 
\end{Proposition}

\proof
Let $U(y) = \{x: |x| \leq y\}$.  So $U=U(e)$.  Let $\widetilde{G} = 
\union_{N \in \Nn} U(Ne)$.
By (1) of \ref{7.8},  $G_\infty \subseteq \widetilde{G}$, and $*,+$ 
coincide on $\widetilde{G}$ up to $H$.  In fact 
$\widetilde{G}=G_\infty$, since $\widetilde{G}/H=\Rr^n$ and $U/H$ 
contains an open neighborhood of 0 in $\Rr^n$.

Since $*,+$ coincide on  $\widetilde{G}$ up to $H$, the proposition 
reduces to the case $G=(R^n,+)$, where $(R,+)$ is the underlying 
additive group of the o-minimal structure.
In this case, add predicates for both
$\{x: (\exists N \in \Nn) |x| \leq N e \}$ and for $\{x: (\forall N 
\in \Nn) |x| < e/N \}$,
obtain weak-o-minimality of their quotient by \cite{BP}, and proceed 
as in the proof of 10.4.
\qed

\begin{Corollary} \label{10.10} $G_\infty^{00} = H$ \end{Corollary}
\proof  A generic set has
generic image in $G_\infty/H$, hence contains a non-small subset of 
$G_\infty/H$, hence the pullback contains at least one full
coset of $H$.
Since $G_\infty/H$ is bounded, $G_\infty^{00}=H$.  \qed

\end{document}